\documentclass[a4paper,oneside,11pt]{article}

\usepackage{amsmath,amsfonts,amscd,amssymb}
\usepackage{longtable,geometry}
\usepackage[english]{babel}
\usepackage[utf8]{inputenc}
\usepackage[active]{srcltx}
\usepackage[T1]{fontenc}
\usepackage{graphicx}
\usepackage{pstricks}
\usepackage{bbm}
\usepackage{tikz}

\usepackage{enumitem}
\usepackage{MnSymbol}
\usepackage{stmaryrd}
\usepackage{nicefrac}
\usepackage{calrsfs}

 \usepackage{rotating}

\usepackage{xcolor}
\usepackage{framed}

\newtheorem{theorem}{Theorem}[section]
\newtheorem{conj}{Conjecture}

\newtheorem{conjecture}[conj]{Conjecture}

\newtheorem{remark}[theorem]{Remark}

\newcommand{\boxRectangle}[3][solid]{\vcenter{\hbox{\;\tikz[scale=1.5]{
    \draw (0, 0) rectangle (6ex, 3ex);
    \draw [#1] (0, 2ex) .. controls (2.6ex, 3ex) and (3.3ex, 1ex) .. (6ex, 2ex);
    \draw (3ex, 1.3ex) node[anchor=north] {$\scriptstyle #2$};
    \draw (5.7ex, 1.5ex) node[anchor=west] {$\scriptstyle #3$};
  }}}
}

\let\OLDthebibliography\thebibliography
\renewcommand\thebibliography[1]{
  \OLDthebibliography{#1}
  \setlength{\parskip}{0pt}
  \setlength{\itemsep}{0pt plus 0.3ex}
}

\title{Sixty years of percolation}
\author{Hugo Duminil-Copin\thanks{
\texttt{duminil@ihes.fr} Institut des Hautes \'Etudes Scientifiques and Universit\'e de Gen\`eve
 \newline
This research was funded by an IDEX Chair from Paris Saclay, by the NCCR SwissMap from the Swiss NSF and the ERC grant 757296 CRIBLAM. We thank David Cimasoni,  S\'ebastien Martineau, Aran Raoufi and Vincent Tassion for their comments on the manuscript.}}
\date{\today}

\begin{document}
\maketitle

\begin{abstract}
Percolation models describe the inside of a porous material. The theory emerged timidly  in the middle of the twentieth century before becoming one of the major objects of interest in probability and mathematical physics. The golden age of percolation is probably the eighties, during which most of the major results were obtained for the most classical of these models, named Bernoulli percolation, but it is really the two following decades which put percolation theory at the crossroad of several domains of mathematics. In this short (and very partial) review, we propose to describe briefly some of the recent progress as well as some famous challenges remaining in the field. 
This review is not intended to probabilists (and a fortiori not to specialists in percolation theory): the target audience is mathematicians of all kinds. \end{abstract}

\section{A brief history of Bernoulli percolation}

\subsection{What is percolation?}

 Intuitively, it is a simplistic probabilistic model for a porous stone. The inside of the stone is described as a random maze in which water can flow. The question then is to understand which part of the stone will be wet when immersed in a bucket of water. Mathematically, the material is modeled as a random subgraph of a reference graph $\mathbb G$ with (countable) vertex-set $\mathbb V$ and edge-set $\mathbb E$ (this is a subset of unordered pairs of elements in $\mathbb V$). 

Percolation on $\mathbb G$ comes in two kinds, {\em bond} or {\em site}. In the former, each edge $e\in \mathbb E$ is either {\em open} or {\em closed}, a fact which is encoded by a function $\omega$ from the set of edges to $\{0,1\}$, where $\omega(e)$ is equal to 1 if the edge $e$ is open, and $0$ if it is closed. We think of an open edge as being open to the passage of water, while closed edges are not. 
A bond percolation model then consists in choosing edges of $\mathbb G$ to be open or closed at random. Site percolation is the same as bond percolation except that, this time,  vertices $v\in \mathbb V$ are either open or closed, and therefore $\omega$ is a (random) function from $\mathbb V$ to $\{0,1\}$.  

The simplest and oldest model of bond percolation, called {\em Bernoulli percolation}, was introduced by Broadbent and Hammersley \cite{BroHam57}. In this model, each edge is open with probability $p$ in $[0,1]$ and therefore closed with probability $1-p$, independently of the state of other edges. Equivalently, the $\omega(e)$ for $e\in\mathbb E$ are independent Bernoulli random variables of parameter $p$.

Probabilists are interested in connectivity properties of the random object obtained by taking the graph induced by $\omega$. In the case of bond percolation, the vertices of this graph are the vertices of $\mathbb G$, and the edges are given by the open edges only. In the case of site percolation, the graph is the subgraph of $\mathbb G$ induced by the open vertices, i.e.~the graph composed of open vertices and edges between them.

%Most early works on Bernoulli percolation focused on the two-dimensional case, where the graph $\mathbb G$ is chosen to be the square lattice $\mathbb Z^2=(\mathbb V^2,\mathbb E^2)$, where $\mathbb V^2$ is the set of points of $\mathbb R^2$ with integer coordinates, and $\mathbb E^2$ is the set of unordered pairs of vertices at Euclidean distance one of each others. 

Let us focus for a moment on Bernoulli percolation on the hypercubic lattice $\mathbb Z^d$ with vertex-set given by the points of $\mathbb R^d$ with integer coordinates, and edges between vertices at Euclidean distance 1 of each other. The simplest connectivity property to study is the fact that the connected component of the origin is finite or not. Set $\theta(p)$  for the probability that the origin is in an infinite connected component of $\omega$. The union bound easily implies that the probability that the origin is connected to distance $n$ is smaller than $(2dp)^n$ (simply use the fact that one of the less than $(2d)^n$ self-avoiding paths of length $n$ starting from the origin must be made of open edges only, as well as the union bound). As a consequence, one deduces that $\theta(p)=0$ as soon as $p<1/(2d)$. This elementary argument was described in the first paper \cite{BroHam57} on percolation theory. On $\mathbb Z^2$, Harris drastically improved this result in \cite{Har60} by showing that $\theta(\tfrac12)=0$.

A slightly harder argument \cite{Ham59} involving Peierls's argument (left to the reader) shows that when $d\ge2$ and $p$ is close to $1$, then $\theta(p)$ is strictly positive. This suggests the existence of a {\em phase transition} in the model: for some values of $p$, connected components are all finite, while for others, there exists an infinite connected component in $\omega$. One can in fact state a more precise result \cite{BroHam57}. For Bernoulli percolation on transitive\footnote{A graph is {\em transitive} if its group of automorphisms acts transitively on its vertices.} (infinite) graphs, there exists $p_c(\mathbb G)\in[0,1]$ such that the probability that there is an infinite connected component in $\omega$ is zero if $p<p_c(\mathbb G)$, and one if $p>p_c(\mathbb G)$ (note that nothing is said about what happens at criticality).
This is an archetypical example of a phase transition in statistical physics: as the parameter $p$ (which can be interpreted physically as the porosity of the stone) is varied continuously through the value $p_c(\mathbb G)$, the probability of  having an infinite connected component jumps from 0 to 1.
 \begin{figure}[h]
\begin{center}\includegraphics[width=0.30\textwidth]{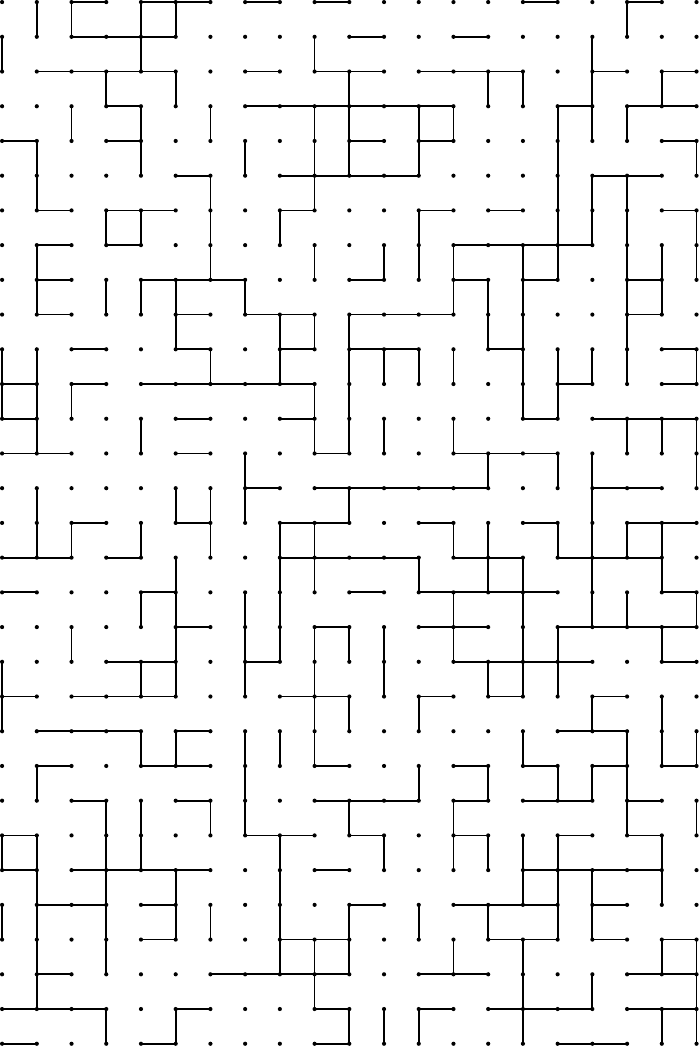}\qquad\includegraphics[width=0.30\textwidth]{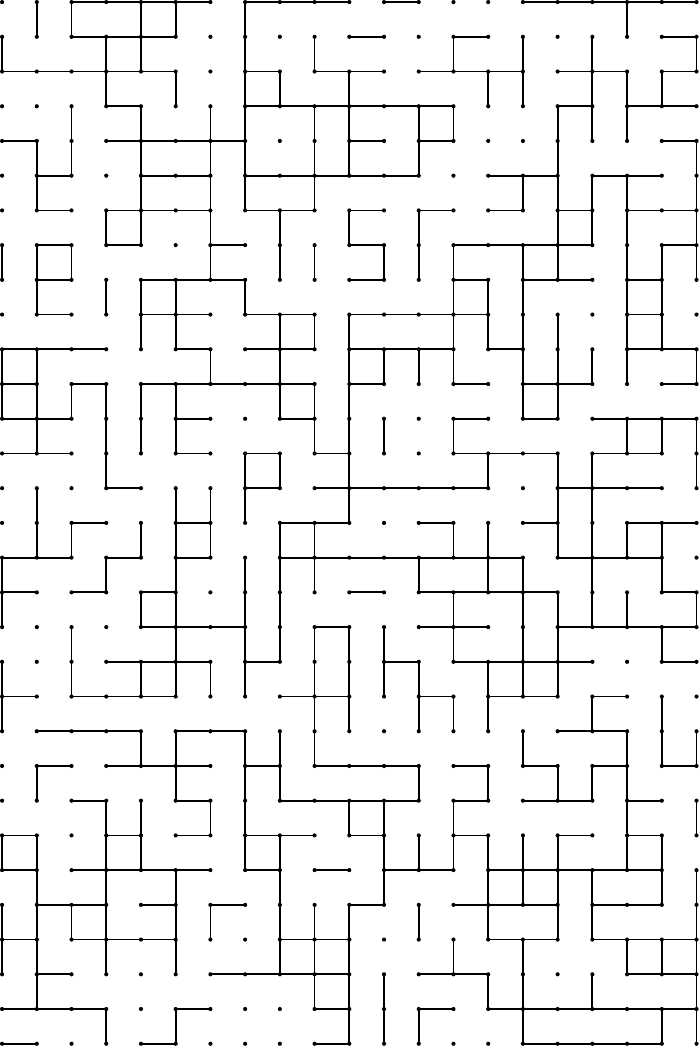}\qquad\includegraphics[width=0.30\textwidth]{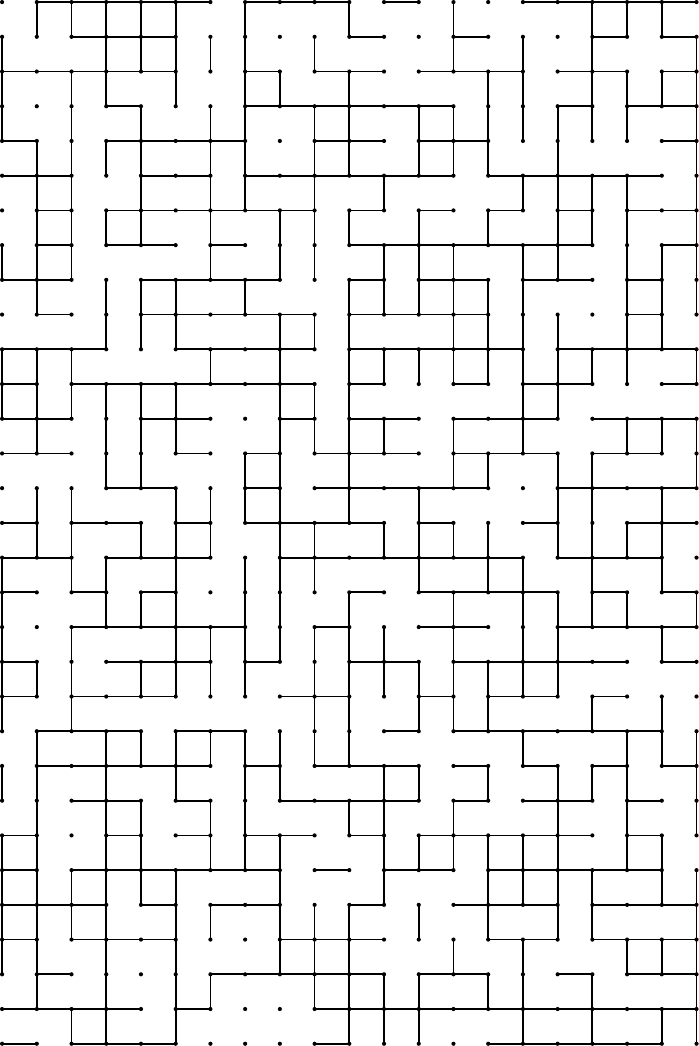}
\caption{A sampled configuration $\omega$ of Bernoulli bond percolation on the square lattice $\mathbb Z^2$ for the values of the parameter $p<1/2$, $p=1/2$ and $p>1/2$.}\end{center}
\end{figure}
\subsection{The eighties: the Golden Age of Bernoulli percolation}

The eighties are famous for pop stars like Michael Jackson and Madonna, and a little bit less for probabilists such as Harry Kesten and Michael Aizenman. Nonetheless, these mathematicians participated intensively in the amazing progress that the theory underwent during this period. 

The decade started in a firework with Kesten's Theorem \cite{Kes80} showing that the critical point of Bernoulli bond percolation on the square lattice $\mathbb Z^2$ is equal to $1/2$.
This problem drove most of the efforts in the field until its final solution by Kesten, and some of the ideas developed in the proof became instrumental  in the thirty years that followed. The strategy is based on an important result obtained simultaneously by Russo \cite{Rus78} and Seymour-Welsh \cite{SeyWel78}, which will be discussed in more details in the next sections. 

While the two-dimensional case concentrated most of the early focus, the model is, of course, not restricted to the graph $\mathbb Z^2$. On $\mathbb Z^d$, Menshikov \cite{Men86} at the same time as Aizenman and Barsky \cite{AizBar87} showed that not only the probability of being connected to distance $n$ is going to 0 when $p<p_c(\mathbb Z^d)$, but that in fact this quantity is decaying exponentially fast, in the sense that there exists $c>0$ depending on $p$ (and $d$) such that $$\theta_n(p)\stackrel{\rm def}=\mathbb P_p[0\text{ is connected to distance $n$}]\le \exp(-cn)$$ for every $n\ge1$.
This result, known under the name of {\em sharpness of the phase transition}, provides a very strong control of the size of connected components. In particular it says that, when $p<p_c$, the largest connected component in a box of size $n$ is typically of size $\log n$. It is the cornerstone of the understanding of percolation in the {\em subcritical regime} $p<p_c$, and as such, represents an important breakthrough in the field. 

Properties of the {\em supercritical regime} $p>p_c(\mathbb Z^d)$ were also studied in detail during this period. For instance, it is natural to ask oneself whether the infinite connected component is unique or not in this regime. In 1987, Aizenman, Kesten and Newman \cite{AizKesNew87} showed that this is indeed the case\footnote{Picturally, in two dimensions,  the infinite connected component has properties similar to those of $\mathbb Z^d$ and can be seen as an ocean. The finite connected components can then be seen as small lakes separated from the sea by the closed edges (which somehow can be seen as the land forming finite islands).}. The proof relied on delicate properties of Bernoulli percolation, and did not extend easily to more general models. Two years later, Burton and Keane proposed a beautiful argument \cite{BurKea89}, which probably deserves its place in The Book, showing by ergodic means that a large class of percolation models has a unique infinite connected component in the supercritical regime. 
A consequence of this theorem is the continuity of $p\mapsto\theta(p)$ when $p>p_c(\mathbb Z^d)$. Of course, many other impressive results concerning the supercritical regime were obtained around the same period, but the lack of space refrains us from describing them in detail.

%Once uniqueness of the infinite connected component has been established, one may ask how fast is the probability that two points are in the same finite connected component decaying. Chayes, Chayes, Grimmett, Newman and Schonmann \cite{ChaChaGriKesSch89} showed that this decay is exponential. This illustrates the fact non-critical percolation is fairly well understood.

The understanding of percolation at $p=p_c(\mathbb Z^d)$ also progressed in the late eighties and in the beginning of the nineties. Combined with the early work of Harris \cite{Har60} who proved $\theta(1/2)=0$ on $\mathbb Z^2$, Kesten's result directly implies that $\theta(p_c)=0$. In dimension $d\ge19$, Hara and Slade \cite{HarSla90} used a technique known under the name of {\em lace expansion} to show that  critical percolation exhibits a {\em mean-field behavior}, meaning that the critical exponents describing the phase transition are matching those predicted by the so-called mean-field approximation. In particular, the mean-field behavior implies that $\theta(p_c)$ is equal to 0. Each few years, more delicate uses of the lace-expansion enable to reduce the dimension starting at which the mean-field behavior can be proved: the best know result today is $d\ge11$ \cite{FitHof15}.

One may wonder whether it would be possible to use the lace expansion to prove that $\theta(p_c)$ is equal to 0 for every dimension $d\ge3$. Interestingly, the mean-field behavior is expected to hold only when $d\ge6$, and to fail for dimensions $d\le 5$ (making the lace expansion obsolete). This leaves the intermediate dimensions $3$, $4$ and $5$ as a beautiful challenge to mathematicians. In particular, the following question is widely considered as the major open question in the field.

\begin{conjecture}\label{conj:1}
Show that $\theta(p_c)=0$ on $\mathbb Z^d$ for every $d\ge3$.
\end{conjecture}
This conjecture, often referred to as the ``$\theta(p_c)=0$ conjecture'', is one of the problems that Harry Kesten was describing in the following terms in his famous 1982 book \cite{Kes82}:
\medbreak
{\em `` Quite apart from the fact that percolation theory had its origin
in an honest applied problem, it is
a source of fascinating problems of the best kind a mathematician can
wish for: problems which are easy to state with a minimum of preparation,
but whose solutions are (apparently) difficult and require new methods. ''}
\medbreak
 It would be unfair to say that the understanding of critical percolation is non-existent for $d\in\{3,4,5\}$. Barsky, Grimmett and Newman \cite{BarGriNew91} proved  that the probability that there exists an infinite connected component in $\mathbb N\times\mathbb Z^{d-1}$ is zero for $p=p_c(\mathbb Z^d)$. It seems like a small step to bootstrap this result to the non-existence of an infinite connected component in the full space $\mathbb Z^d$... But it is not. More than twenty five years after \cite{BarGriNew91}, the conjecture still resists  and basically no improvement has been obtained.

\subsection{The nineties: the emergence of new techniques}\label{sec:tool}

Percolation theory underwent a major mutation in the 90's and early 00's. While some of the historical questions were solved in the previous decade, new challenges appeared at the periphery of the theory. In particular, it became clear that a deeper understanding of percolation would require the use of techniques coming from a much larger range of mathematics. As a consequence, Bernoulli percolation took a new place at the crossroad of several domains of mathematics, a place in which it is not just a probabilistic model anymore.

\subsubsection{Percolation on groups} In a beautiful paper entitled {\em Percolation beyond $\mathbb Z^d$, many questions and a few answers} \cite{BenSch96}, Benjamini and Schramm underlined the relevance of Bernoulli percolation on Cayley graphs\footnote{The Cayley graph $\mathbb G=\mathbb G(G,S)$ of a finitely generated group $G$ with a symmetric system of generators $S$ is the graph with vertex-set $\mathbb V=G$ and edge-set given by the unordered pairs $\{x,y\}\subset G$ such that $yx^{-1}\in S$. For instance, $\mathbb Z^d$ is a Cayley graph for the free abelian group with $d$ generators.} of finitely generated infinite groups by proposing a list of problems relating properties of Bernoulli percolation to properties of groups. The paper triggered a number of new problems in the field and drew the attention of the community working in geometric group theory on the potential applications of percolation theory. 
%The paper proposed a list of problem relating the properties of Bernoulli percolation to those of the group $G$. 

A striking example of a connection between the behavior of percolation and properties of groups is  provided by the following conjecture, known as the ``$p_c<p_u$ conjecture''. Let $p_u(\mathbb G)$ be the smallest value of $p$ for which the probability that there exists a {\em unique} infinite connected component is one. On the one hand, the uniqueness result \cite{BurKea89} mentioned in the previous section implies that $p_c(\mathbb Z^d)=p_u(\mathbb Z^d)$. On the other hand, one can easily convince oneself that, on an infinite tree $\mathbb T_d$ in which each vertex is of degree $d+1$, one has $p_c(\mathbb T_d)=1/d$ and $p_u(\mathbb T_d)=1$. More generally, the $p_c<p_u$ conjecture relates the possibility of infinitely many connected components to the property of non-amenability\footnote{$G$ is {\em non-amenable} if for any Cayley graph $\mathbb G$ of $G$, the infimum of $|\partial A|/|A|$ on non-empty finite subsets $A$ of $G$ is strictly positive, where $\partial A$ denotes the {\em boundary} of $A$ (i.e.~the set of $x\in A$ having one neighbor outside $A$)  and $|B|$ is the cardinality of the set $B$. }
 of the underlying group.\begin{conjecture}[Benjamini-Schramm]
Consider a Cayley graph $\mathbb G$ of a finitely generated (infinite) group $G$. Then 
$$p_c(\mathbb G)<p_u(\mathbb G)\ \Longleftrightarrow\text{ $G$ is non-amenable}.$$
\end{conjecture}

The most impressive progress towards this conjecture was achieved by Pak and Smirnova \cite{PakSmi00} who provided a group theoretical argument showing that any non-amenable group possesses {\em a} multi-system of generators for which the corresponding Cayley graph satisfies $p_c<p_u$. This is a perfect example of an application of Geometric Group Theory to probability. Nicely enough, the story sometimes goes the other way and percolation shed a new light on classical questions on group theory. The following example perfectly illustrates this cross-fertilization. 

In 2009, Gaboriau and Lyons \cite{GabLyo09} provided a measurable solution to von Neumann's (and Day's)  famous problem on non-amenable groups. %The notion of non-amenable group was introduced by von Neumann to explain the Banach-Tarsky paradox. 
While it is simple to show that a group containing the free group $F_2$ as a subgroup is non-amenable, it is non-trivial to determine whether the converse is true. Olsanski \cite{Ols80} showed in 1980 that this is not the case, but Whyte \cite{Why99} gave a very satisfactory geometric solution: a
finitely generated group is non-amenable if and only if it admits a partition into pieces that are all uniformly bi-lipschitz
equivalent to the regular four-valent tree $\mathbb T_3$. Bernoulli percolation was used by Gaboriau and Lyons to show the measurable counterpart of this theorem, a result which has many important applications in the ergodic theory of group actions. %Other examples of applications of Bernoulli percolation to geometric group theory can be found in \cite{??,??,??}. 

\begin{figure}
\begin{center}\includegraphics[width=0.45\textwidth]{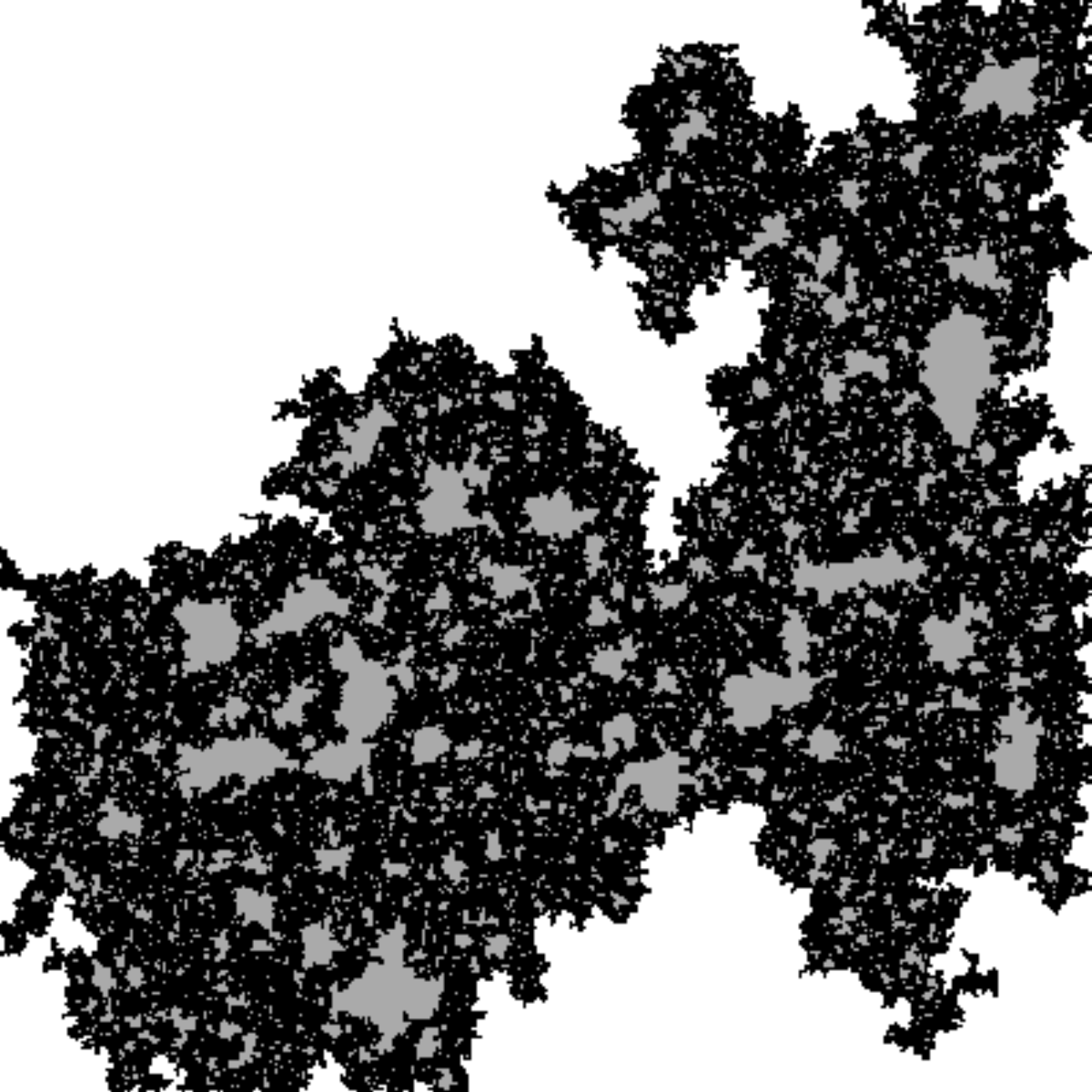}\caption{\label{fig:algo} A large connected component of $\omega$ for critical bond percolation on the square lattice $\mathbb Z^2$ (simulation by Vincent Beffara).}\end{center}
\end{figure}

\subsubsection{Complex analysis and percolation}\label{sec:11}
The nineties also saw a renewed interest in questions about planar percolation. The impressive developments in the eighties of Conformal Field Theory, initiated by  Belavin, Polyakov and Zamolodchikov \cite{BelPolZam84},  suggested that the scaling limit of planar Bernoulli percolation is conformally invariant at criticality. From a mathematical perspective, the notion of conformal invariance of the
entire model is ill-posed, since the meaning of scaling limit depends on the object
under study (interfaces, size of connected components, crossings, etc). In 1992, the observation that properties
of interfaces should also be conformally invariant led Langlands, Pouliot and
Saint-Aubin \cite{LanPouSai94} to publish numerical values in agreement with the conformal
invariance in the scaling limit of crossing probabilities; see Section~\ref{sec:3} for a precise definition (the authors attribute this conjecture to  Aizenman). The same year, the physicist Cardy \cite{Car92}
proposed an explicit formula for the limit of crossing probabilities in rectangles of fixed aspect ratio. 

These two papers, while numerical (for the first one) and physical (for the second one), attracted many mathematicians
to the domain.  
In 2001, Smirnov \cite{Smi01} proved  Cardy's
formula for critical site percolation on the triangular lattice, hence
rigorously providing a concrete example of a conformally invariant property of
the model. In parallel, a major breakthrough shook the field of probability. In 2000, Schramm \cite{Sch00} introduced a random process, now called the {\em Schramm-Loewner Evolution}, describing the behavior of interfaces between open and closed sites. Very shortly after Smirnov's proof of Cardy's formula, the complete description of the scaling limit of site percolation was obtained, including the description of the full ``loop ensemble'' corresponding to the interfaces bordering each connected component by Camia and Newman \cite{CamNew06}, see \cite{BefDum13} for more references on this beautiful subject. 

Smirnov's theorem and the Schramm-Loewner Evolution share a common feature: they both rely on complex analysis in a deep way. The first result uses discrete complex analysis, i.e.~the study of functions on graphs that approximate holomorphic functions, to prove the convergence of certain observables of the model to conformal maps\footnote{Nicely enough, the story can again go the other way: Smirnov's argument can also be used to provide an alternative proof of Riemann's mapping theorem \cite{Smi17}.
}. The second revisits Loewner's deterministic evolutions (which were used to solve the Bieberbach conjecture) to construct random processes whose applications now spread over all probability theory. 
\subsubsection{Discrete Fourier analysis, concentration inequalities... and percolation} The end of the nineties witnessed the appearance of two important new problems regarding Bernoulli percolation. H\"aggstr\"om, Peres and Steif \cite{HagPerSte97} introduced a simple time dynamics in Bernoulli percolation obtained by resampling each edge at an exponential rate 1. More precisely, an exponential clock is attached to each edge of the lattice, and each time the clock of an edge rings, the state of the edge is resampled.  It turns out that dynamical percolation is a very interesting model which exhibits quite rich phenomena. 

As mentioned above, it is known since Harris that $\theta(1/2)=0$ on $\mathbb Z^2$. Since Bernoulli percolation is the invariant measure for the dynamics, this easily implies that for any fixed time $t\ge0$, the probability that dynamical percolation at time $t$ does not contain an infinite connected component is zero. Fubini's theorem shows a stronger result: with probability 1, the set of times at which the configuration contains an infinite connected component is of Lebesgue measure 0. Nonetheless, this statement does not exclude the possibility that this set of times is non-empty.%Proving the existence of such exceptional times has been a question of major interest in percolation theory. 

In 1999, 
Benjamini, Kalai, and Schramm \cite{BenKalSch99} initiated the study of the Fourier spectrum of percolation and its applications to the noise sensitivity of the model measuring how much connectivity properties of the model are robust under the dynamics. This work advertised the usefulness of concentration inequalities and discrete Fourier analysis for the understanding of percolation: they provide information on the model which is invisible from the historical probabilistic and geometric approaches. We will see below that these tools will be crucial to the developments of percolation theory of dependent models.  

We cannot conclude this section without mentioning the impressive body of works by Garban, Pete and Schramm \cite{GarPetSch11,GarPetSch13,GarPetSch17} describing in detail the noise sensitivity and dynamical properties of planar percolation. These works combine the finest results on planar Bernoulli percolation and provide a precise description of the behavior of the model, in particular proving that the Hausdorff dimension of the set of times at which there exists an infinite connected component is equal to $31/36$. 
\bigbreak
{\em At this stage of the review, we hope that the reader gathered some understanding of the problematic about Bernoulli percolation, and got some idea of the variety of fields of mathematics involved in the study of the model. We will now try to motivate the introduction of more complicated percolation models before explaining, in Section~\ref{sec:4}, how some of the techniques mentioned above can be used to study these more general models. We refer to \cite{Gri99} for more references.}

\section{Beyond Bernoulli percolation}
While the theory of Bernoulli percolation still contains a few gems that deserve a solution worthy of their beauty, recent years have revived the interest for more general percolation models appearing in various areas of statistical physics as natural models associated with other random systems. While Bernoulli percolation is a product measure, the states of edges in these percolation models are typically not independent. 
Let us discuss a few ways of introducing these ``dependent'' percolation models.

\subsection{From spin models to bond percolation}

Dependent percolation models are often associated with lattice spin models. These models have been introduced as discrete models for real life experiments and have been later on found useful to model a large variety of phenomena and systems ranging from ferromagnetic materials to lattice gas. They also provide discretizations of Euclidean and Quantum Field Theories and are as such important from the point of view of theoretical physics. While the original motivation came from physics, they appeared as extremely complex and rich mathematical objects, whose study required the developments of new tools that found applications in many other domains of mathematics. 

The archetypical example of the relation spin model/percolation is provided by the {\em Potts model} and FK percolation (defined below). In the former, spins are chosen among a set of $q$ colors (when $q=2$, the model is called the {\em Ising model}, and spins are seen as $\pm1$ variables), and the distribution depends on a parameter $\beta$ called the {\em inverse temperature}. We prefer not to define the models too formally here and refer to \cite{Dum17} for more details. We rather focus on the relation between these models and a dependent {\em bond} percolation model called {\em Fortuin-Kasteleyn (FK) percolation}, or {\em random-cluster model}.

FK percolation was
 introduced by Fortuin and Kasteleyn in \cite{ForKas72} as a unification of different models of statistical physics satisfying series/parallel laws when modifying the underlying graph. The probability measure on a finite graph $\mathbb G$, denoted by 
$\mathbb P_{\mathbb G,p,q}$, depends on two parameters -- the {\em edge-weight} 
$p\in[0,1]$ and the {\em cluster-weight} $q>0$ -- and is defined by the formula
\begin{equation}
  \label{probconf}
  \mathbb P_{\mathbb G,p,q}[\{\omega\}] :=
  \frac {p^{|\omega|}(1-p)^{|\mathbb E|-|\omega|}q^{k(\omega)}}
  {Z(\mathbb G,p,q)}\qquad\text{ for every }\omega\in\{0,1\}^\mathbb E,
\end{equation}
where $|\omega|$ is the number of open edges and $k(\omega)$ the number of connected components in $\omega$. The constant $Z(\mathbb G,p,q)$ is a 
normalizing constant, referred to as the \emph{partition function}, defined in such a way that the sum over all configurations equals 1. 
For $q=1$, the model is Bernoulli percolation, but for $q\ne1$, the probability measure is different, due to the term $q^{k(\omega)}$ taking into account the number of connected components. 
Note that, at first sight, the definition of the model makes no sense on infinite graphs (contrarily to Bernoulli percolation) since the number of open edges (or infinite connected components) would be infinite. Nonetheless, the model can be extended to infinite graphs by taking the (weak) limit of measures on finite graphs. 
For more on the model, we refer to the comprehensive reference book \cite{Gri06}.

As mentioned above, FK percolation is connected to the Ising and the Potts models via what is known as the {\em Edwards-Sokal coupling}. It is straightforward to describe this coupling in words. Let $\omega$ be a percolation configuration sampled according to the FK percolation with edge-weight $p\in[0,1]$ and cluster-weight $q\in\{2,3,4,\dots\}$. The random coloring of $\mathbb V$ obtained by assigning to connected components of $\omega$ a color chosen uniformly (among $q$ fixed colors) and independently for each connected component, and then giving to a vertex the color of its connected component, is a realization of the $q$-state Potts model at inverse temperature $\beta=-\tfrac12\log(1-p)$. For $q=2$, the colors can be understood as $-1$ and $+1$ and one ends up with the Ising model.
 
The relation between FK percolation and the Potts model is not an exception. Many other lattice spin models also possess their own Edwards-Sokal coupling with a dependent percolation model. This provides us with a whole family of natural dependent percolation models that are particularly interesting to study. We refer the reader to \cite{ChaMac97,PfiVel97} for more examples.
%Below, we wish to highlight a new example of such a relation, which motivates the introduction of another percolation model which will come back naturally in the next sections. 

%Define the {\em XOR Ising model} obtained by considering the spin configuration $(\sigma_x\tau_x)_{x\in \mathbb V}$, where $\sigma$ and $\tau$ are two independent Ising configurations at inverse temperature $\beta$. In the XOR Ising model, the spin (at $x$) takes value $+1$ if the spins at $x$ are equal in $\sigma$ and $\tau$, and spin $-1$ if they are different. The name XOR Ising model was coined recently by Wilson \cite{Wil11}. Now, define the {\em double random current model} to be the percolation model obtained as follows: the probability of $\omega\in\{0,1\}^\mathbb E$ is $$\mathbb P[\{\mathbf \omega\}]=\frac{p^{|\omega|}(1-p)^{|\mathbb E|-|\omega|}\cdot2^{k(\omega)+|\omega|}\sum_{\eta\subset \omega} p^{-|\eta|/2}}{Z(\mathbb G,p)},$$where the sum is over even subgraphs $\eta$, i.e.~subgraphs having even degree at every vertex. The XOR Ising model can be recovered from the double random current model by choosing a spin $+1$ or $-1$ uniformly for every connected component of $\omega$, exactly like the Ising model is obtained from FK percolation with cluster-weight $2$. Note that the structure of the double random current model is the same as for FK percolation: there is an energy term of the form $p^{|\omega|}(1-p)^{|\mathbb E|-|\omega|}$, and a second term taking into account the global connectivity properties of the configuration.

\subsection{From loop models to site percolation models}\label{sec:2}

In two dimensions, there is another recipe to obtain dependent percolation models, but on sites this time. %Note that the Ising model or the XOR Ising model are examples of site percolation models.
Each site percolation configuration on a {\em planar graph} is naturally associated, by the so-called {\em low-temperature expansion}, with a bond percolation of a very special kind, called a {\em loop model}, on the dual graph $\mathbb G^*=(\mathbb V^*,\mathbb E^*)$, where $\mathbb V^*$ is the set of faces of $\mathbb G$, and $\mathbb E^*$ the set of unordered pairs of adjacent faces. More precisely, if $\omega$ is an element of $\{0,1\}^\mathbb V$ (which can be seen as an attribution of $0$ or $1$ to faces of $\mathbb G^*$), define the percolation configuration $\eta\in\{0,1\}^{\mathbb E^*}$ by first extending $\omega$ to the exterior face $x$ of $\mathbb G^*$ by arbitrarily setting $\omega_x=1$, and then saying that an edge of $\mathbb G^*$ is open if the two faces $x$ and $y$ of $\mathbb G^*$ that it borders satisfy $\omega_x\ne\omega_y$. In physics terminology, the configuration $\eta$ corresponds to the {\em domain walls} of $\omega$. 
Notice that the degree of $\eta$ at every vertex is necessarily even, and that therefore $\eta$ is necessarily an even subgraph. 

But one may go the other way: to any even subgraph of $\mathbb G^*$, one may associate a percolation configuration on $\mathbb V$ by setting $+1$ for the exterior face of $\mathbb G^*$, and then attributing 0/1 values to the other faces of $\mathbb G^*$ by switching between 0 and 1  when crossing an edge. This reverse procedure provides us with a recipe to construct new dependent site percolation models: construct first a loop model, and then look at the percolation model it creates.

\begin{figure}
\includegraphics[width=0.27\textwidth,angle=90]{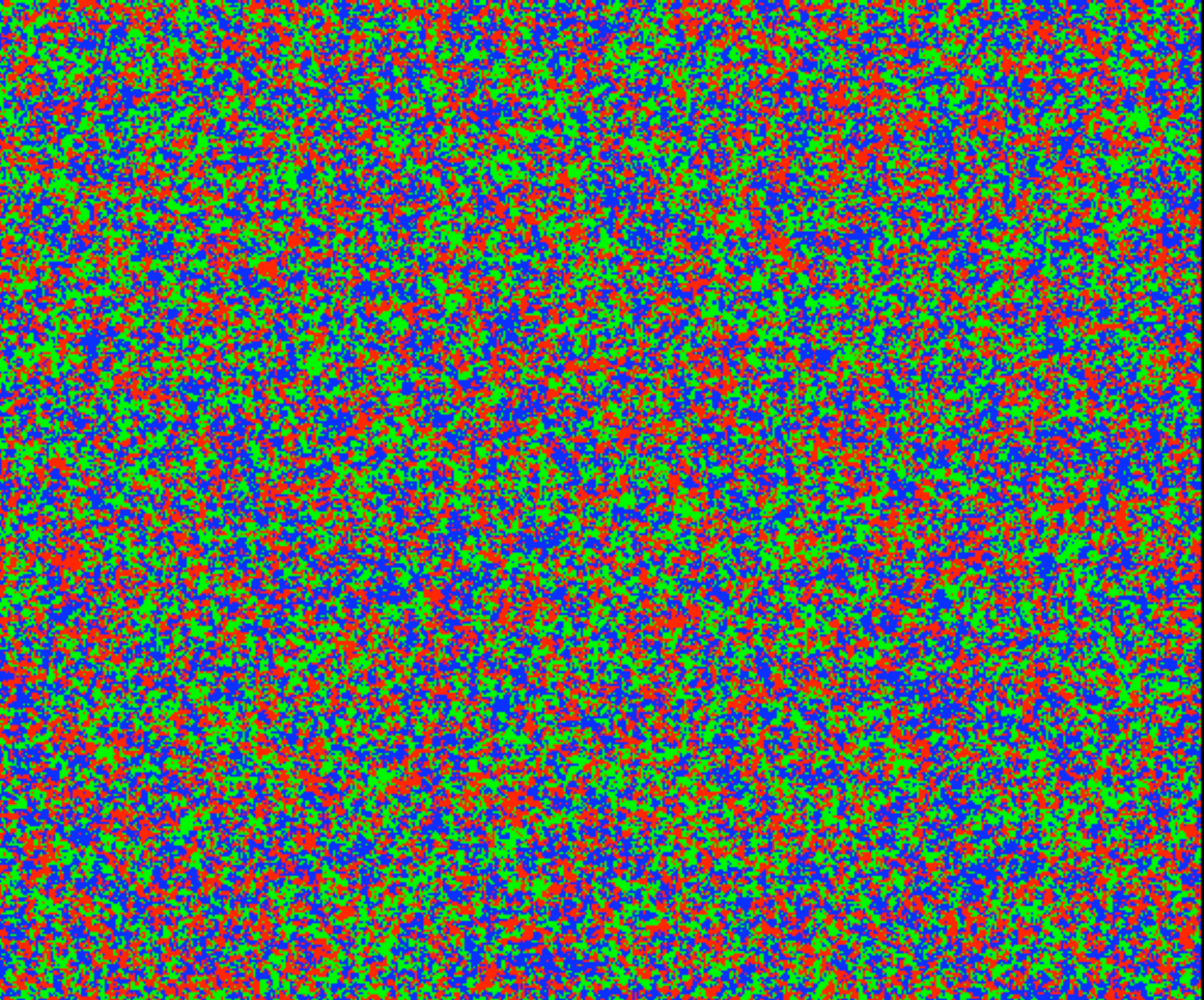}
\includegraphics[width=0.27\textwidth,angle=90]{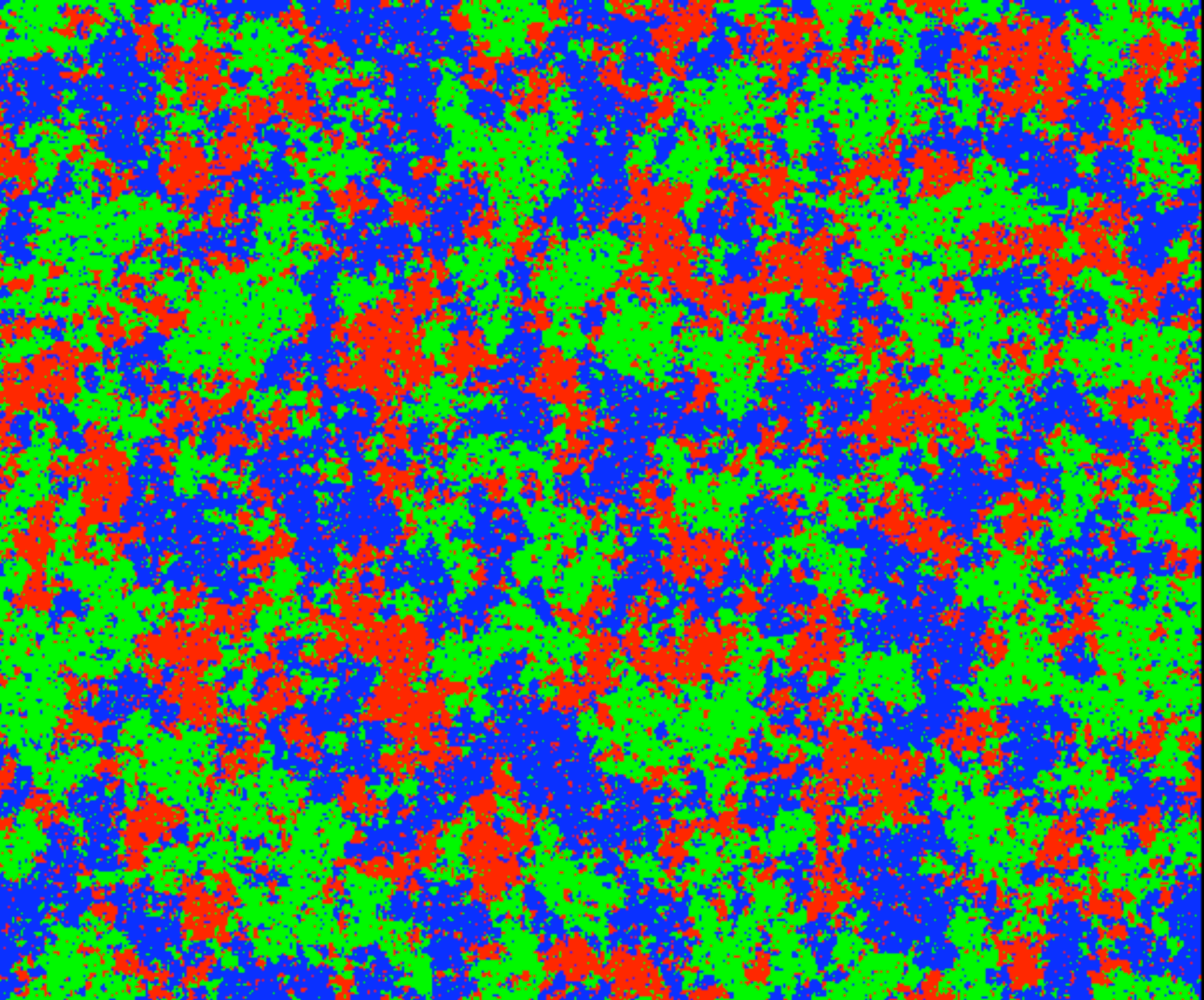}
\includegraphics[width=0.27\textwidth,angle=90]{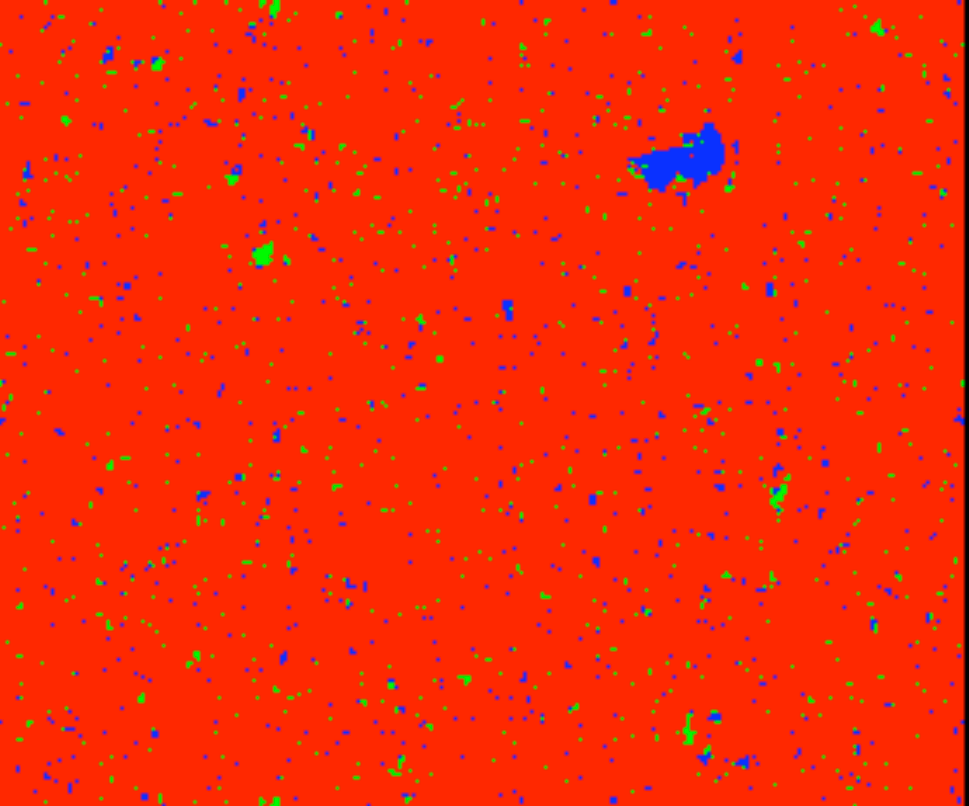}
\caption{Simulations by Vincent Beffara of the three-state planar Potts model obtained from the FK percolation with parameter $p<p_c$, $p=p_c$ and $p>p_c$. On the right, every vertex of the infinite connected component receives the same color, therefore one of the colors wins over the other ones, while this is not the case at criticality or below it.}
\end{figure}

When starting with the Ising model on the triangular lattice (which is indeed a site percolation model: a vertex is open if the spin is $+1$, and closed if it is $-1$), the low-temperature expansion gives rise to a model of random loops on the hexagonal lattice, for which the weight of an even subgraph $\eta$ is proportional to $\exp(-2\beta|\eta|)$.
This loop model  was generalized by Nienhuis et al \cite{DomMukNie81}  to give the {\em loop $O(n)$ model} depending on two parameters, an {\em edge-weight} $x>0$ and a {\em loop-weight} $n\ge0$. It is defined on the hexagonal lattice $\mathbb H=(\mathbb V,\mathbb E)$ as follows: the probability of $\eta$ on $\mathbb H$ is given by
$$\mu_{\mathbb G,x,n}[\{\eta\}]=\frac{x^{|\eta|}n^{\ell(\eta)}}{Z(\mathbb G,x,n)}$$
(where $\ell(\eta)$ is the number of loops in $\eta$) if $\eta$ is an even subgraph, and 0 otherwise.

From this loop model, one obtains a site percolation model on the triangular lattice resembling FK percolation. We will call this model the {\em FK representation of the dilute Potts model} (for $n=1$, it is simply the Ising model mentioned above).
\medbreak
{\em We hope that this section underlined the relevance of {\em some} dependent percolation models, and that the previous one motivated questions for Bernoulli percolation that possess natural counterparts for these dependent percolation models. The next sections describe the developments and solutions to these questions.}

\section{Exponential decay of correlations in the subcritical regime}\label{sec:4}

As mentioned in the first section, proving exponential decay of $\theta_n(p)$ when $p<p_c$ was a milestone in the theory of Bernoulli percolation since it was the key to a deep understanding of the subcritical regime. The goal of this section is to discuss the natural generalizations of these statements to FK percolation with cluster-weight $q>1$. 
Below, we set $\theta_n(p,q)$ and $\theta(p,q)$ for the probabilities of being connected to distance $n$ and to infinity respectively. Also, we define 
\begin{align*}p_c(q)&\stackrel{\rm def}=\inf\{p\in[0,1]:\theta(p,q)>0\},\\
p_{\rm exp}(q)&\stackrel{\rm def}=\sup\{p\in[0,1]:\exists c>0,\forall n\ge0,\,\theta_n(p,q)\le\exp(-c n)\}.\end{align*}
Exponential decay in the subcritical regime gets rephrased as $p_c(q)=p_{\rm exp}(q)$. Exactly like in the case of Bernoulli percolation, the result was first proved in two dimensions, and then in higher dimensions, so that we start by the former. Interestingly, both proofs (the two-dimensional one and the higher dimensional one) rely on the analysis of functions on graphs (discrete Fourier analysis and the theory of randomized algorithms respectively).

\subsection{Two dimensions: studying crossing probabilities}\label{sec:3}

\paragraph{Bernoulli percolation}Let us start by discussing Bernoulli percolation to understand what is going on in a fairly simple case. It appeared quickly (for instance in \cite{Har60}) that crossing probabilities play a central role in the study of planar percolation. For future reference, let ${\bf H}_{n,m}(p)$ be the probability that there exists a crossing from left to right of the rectangle $[0,n]\times[0,m]$, i.e.~a path of open edges (of the rectangle) from $\{0\}\times[0,m]$ to $\{n\}\times[0,m]$, which can be written pictorially as follows
$${\bf H}_{n,m}(p)\stackrel{\rm def}{=}\mathbb P_p\Big[ \boxRectangle n m\Big].$$
It is fairly elementary to prove that:
\begin{itemize}
\item there exists $\varepsilon>0$ such that if ${\bf H}_{k,2k}(p)<\varepsilon$ for some $k$, then $\theta_n(p)$ decays exponentially fast. In words, if the probability of crossing {\em some} wide rectangle in the easy direction is very small, then we are at a value of $p$ for which exponential decay occurs. This implies that for any $p> p_{\rm exp}$, ${\bf H}_{k,2k}(p)\ge \varepsilon$ for every $k\ge1$.
\item there exists $\varepsilon>0$ such that if ${\bf H}_{2k,k}(p)>1-\varepsilon$, then $\theta(p)>0$. Again, in words,  if the probability of crossing {\em some} wide rectangle in the hard direction tends to 1, then we are above criticality. This implies that when $p< p_c$, ${\bf H}_{2k,k}(p)\le 1-\varepsilon$ for every $k\ge1$.
\end{itemize}
In order to prove that the phase transition is sharp, one should therefore prove that there cannot be a whole interval of values of $p$ for which crossings in the easy (resp.~hard) direction occur with probability bounded away uniformly from 0 (resp.~1). The proof involves two important ingredients.

The first one is a result due to 
Russo \cite{Rus78} and Seymour-Welsh \cite{SeyWel78}, today known as the RSW theory. The theorem states that crossing probabilities in rectangles with different aspect ratios can be bounded in terms of each other. More precisely, it shows that for any $\varepsilon>0$ and $\rho>0$, there exists $C=C(\varepsilon,\rho)>0$ such that for every $p\in[\varepsilon,1-\varepsilon]$ and $n\ge1$,
$${\bf H}_{n,n}(p)^C\le {\bf H}_{\rho n,n}(p)\le 1-(1-{\bf H}_{n,n}(p))^C.$$
This statement has a direct consequence: if for $p$, probabilities of crossing rectangles in the easy direction are not going to 0, then the same holds for squares. Similarly, if probabilities of crossing rectangles in the hard direction are not going to 1, then the same holds for squares.
At the light of the previous paragraph, what we really want to exclude now  is the existence of a whole interval of values of $p$ for which the probabilities of crossing squares remain bounded away from 0 and 1 uniformly in $n$.

In order to exclude this possibility, we invoke a second ingredient, which is of a very different nature. It consists in proving that probabilities of crossing squares go quickly from a value close to 0 to a value close to 1. Kesten originally proved this result by hand by showing that the derivative of the function $p\mapsto{\bf H}_{n,n}(p)$ satisfies a differential inequality of the form
\begin{equation}\label{eq:diff}{\bf H}_{n,n}'\ge c\log n \cdot {\bf H}_{n,n}(1-{\bf H}_{n,n}).\end{equation}
This differential inequality immediately shows that the plot of the function ${\bf H}_{n,n}$ has an $S$ shape as on Fig.~\ref{fig:algo}, and that ${\bf H}_{n,n}(p)$ therefore goes from $\varepsilon$ to $1-\varepsilon$ in an interval of $p$ of order $O(1/\log n)$. In particular, it implies that only one value of $p$ can be such that ${\bf H}_{n,n}(p)$ remains bounded away from 0 and 1 uniformly in $n$, hence concluding the proof.

In \cite{BolRio06}, Bollob\'as and Riordan proposed an alternative strategy to prove \eqref{eq:diff}. They suggested using a concept long known to combinatorics: a finite random system undergoes a {\em sharp threshold} if its qualitative behavior changes quickly as the result of a small perturbation of the parameters ruling the probabilistic structure. The notion of sharp threshold emerged in the combinatorics community studying graph properties of random graphs, in particular in the work of Erd\"os and Renyi  \cite{ErdRen59} investigating the properties of Bernoulli percolation on the complete graph. 

Historically, the general theory of sharp thresholds for discrete product spaces was developed by Kahn, Kalai and Linial in \cite{KahKalLin88} in the case of the uniform measure on $\{0,1\}^n$, i.e.~in the case of $\mathbb P_p$ with  $p=1/2$ (see also an earlier non-quantitative version by \cite{Rus82}). There, the authors used the Bonami-Beckner inequality \cite{Bek75,Bon70} together with discrete Fourier analysis to deduce inequalities between the variance of a Boolean function and the covariances (often called {\em influences}) of this function with the random variables $\omega(e)$. Bourgain, Kahn, Kalai,
Katznelson and Linial \cite{BouKahKal92} extended these inequalities to product spaces $[0,1]^n$ endowed with the uniform measure (see also Talagrand \cite{Tal94}), a fact which enables to cover the case of $\mathbb P_p$ for every value of $p\in[0,1]$. 

Roughly speaking, the statement can be read as follows: for any increasing\footnote{Here, increasing is meant with respect to the partial order on $\{0,1\}^{\mathbb E}$ defined by $\omega\le \omega'$ if $\omega(e)\le \omega'(e)$ for every edge $e\in\mathbb E$. Then, $f$ is {\em increasing} if $\omega\le \omega'$ implies $f(\omega)\le f(\omega')$.} (boolean) function ${\bf f}:\{0,1\}^\mathbb E\rightarrow\{0,1\}$,
\begin{equation}\label{eq:a1}{\rm Var}_p({\bf f})\le c(p)\sum_{e\in \mathbb E}\frac{{\rm Cov}_p[{\bf f},\omega(e)]}{\log(1/{\rm Cov}_p[{\bf f},\omega(e)])},\end{equation}
where $c$ is an explicit function of $p$ that remains bounded away from 0 and 1 when $p$ is away from 0 and 1. 

Together with the following differential formula (which can be obtained by simply differentiating $\mathbb E_p[{\bf f}]$)
\begin{equation}\label{eq:a2}\frac{\rm d}{{\rm d}p}\mathbb E_p[{\bf f}]=\frac1{p(1-p)}\sum_{e\in\mathbb E}{\rm Cov}_p[{\bf f},\omega(e)]\end{equation}
for the indicator function {\bf f} of the event that the square $[0,n]^2$ is crossed horizontally, we deduce that 
\begin{equation}\label{eq:a3}{\bf H}_{n,n}'\ge \frac{4}{c(p)\log(1/\max_e{\rm Cov}_p[{\bf f}_n,\omega(e)])} \cdot {\bf H}_{n,n}(1-{\bf H}_{n,n}).\end{equation}
This inequality can be used as follows. The covariance between the existence of an open path and an edge $\omega(e)$ can easily be bounded by the fact that one of the two endpoints of the edge $e$ is connected to distance $n/2$ (indeed, for $\omega(e)$ to influence the outcome of ${\bf f}_n$, there must be an open crossing going through $e$ when $e$ is open). But, in the regime where crossing probabilities are bounded away from 1,  the probability of being connected to distance $n/2$ can easily be proved to decay polynomially fast, so that in fact
${\bf H}'_{n,n}\ge c\log n\cdot {\bf H}_{n,n}(1-{\bf H}_{n,n})$ as required. 
\paragraph{FK percolation}
What survives for dependent percolation models such as FK percolation? The good news is that the BKKKL result can be extended to this context \cite{GraGri06}
to obtain \eqref{eq:a1}. Equation \eqref{eq:a2} is obtained in the same way by elementary differentiation. It is therefore the RSW result which is tricky to extend. 

While mathematicians are in possession of many proofs of this theorem for Bernoulli percolation \cite{BolRio06,BolRio06c,BolRio10,Rus78,SeyWel78,Tas14,Tas14b}, one had to wait for thirty years to obtain the first proof of this theorem for dependent percolation model. 

The following result is the most advanced result in this direction. Let ${\bf H}_{n,m}(p,q)$ be the probability that  $[0,n]\times[0,m]$ is crossed horizontally for FK percolation. 
\begin{theorem}
For any $\rho>0$, there exists a constant $C=C(\rho,\varepsilon,q)>0$ such that for every $p\in[\varepsilon,1-\varepsilon]$ and $n\ge1$,
$${\bf H}_{n,n}(p,q)^C\le{\bf H}_{\rho n,n}(p,q)\le1-(1-{\bf H}_{n,n}(p,q))^C.$$
\end{theorem}

A consequence of all this is the following result \cite{BefDum12} (see also  \cite{DumMan16,DumRaoTas16}).

\begin{theorem}
Consider the FK model with cluster weight $q\ge1$. Then, for any $p<p_c$, there exists $c=c(p,q)>0$ such that for every $n\ge1$,
$$\theta_n(p,q)\le \exp(-cn).$$
\end{theorem}

\begin{figure}
\begin{center}\includegraphics[width=0.35\textwidth]{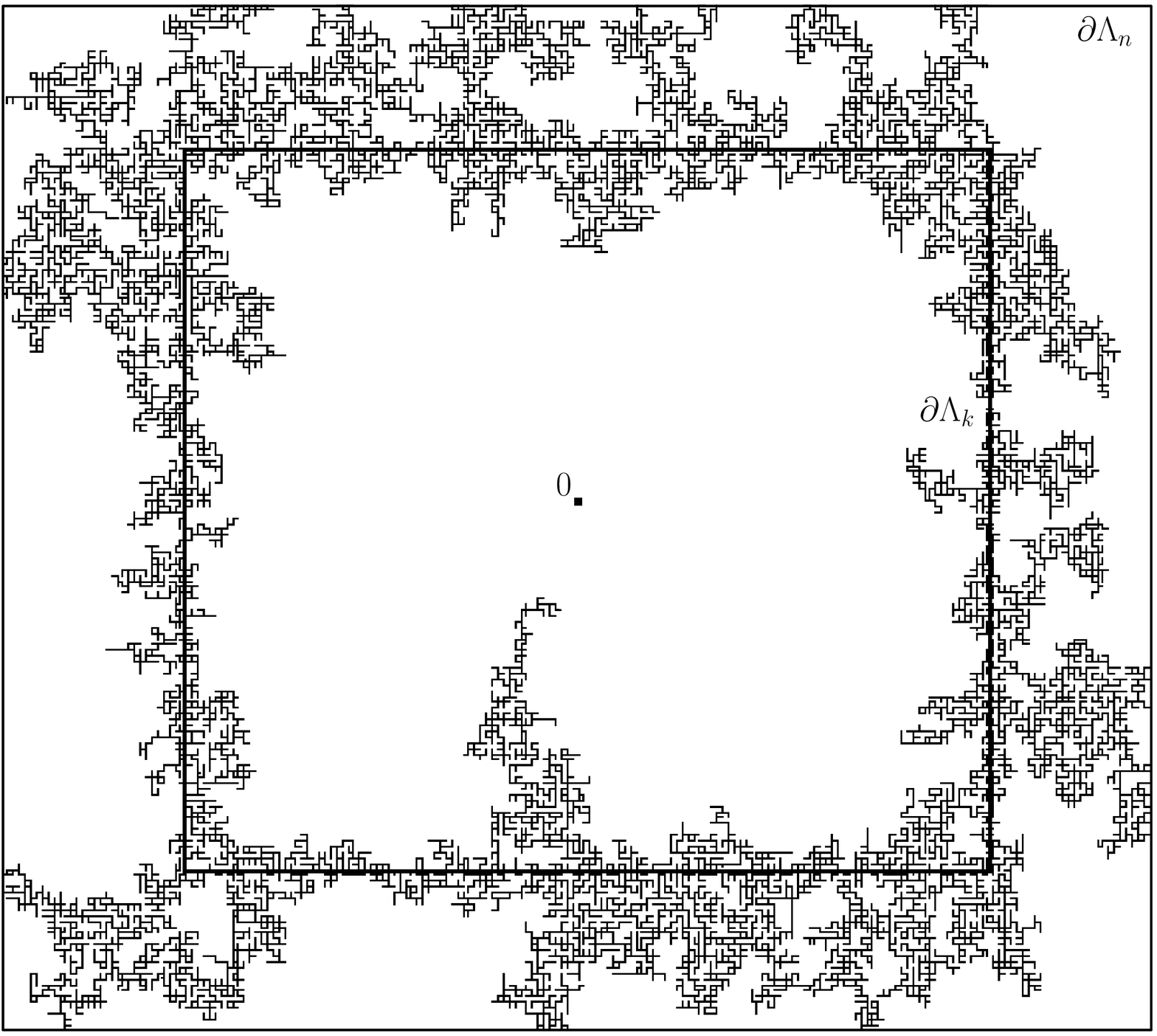}\qquad\qquad\qquad\includegraphics[width=0.30\textwidth]{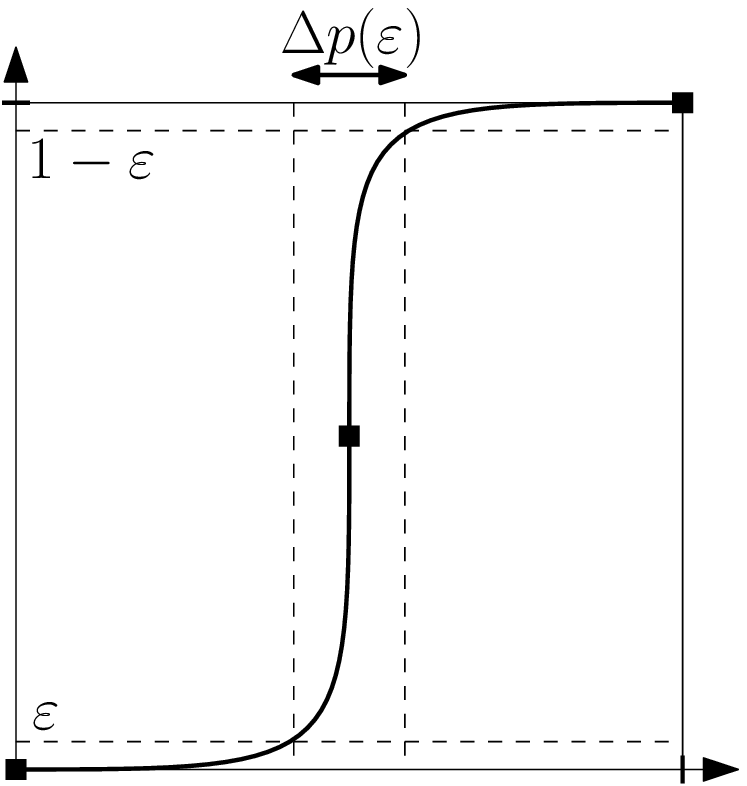}
\caption{\label{fig:algo} {\bf Left.} The randomized algorithm is obtained as follows: pick a distance $k$ to the origin uniformly and then explore ``from the inside'' all the connected components intersecting the boundary of the box of size $k$. If one of these connected components intersect both 0 and the boundary of the box of size $n$, then we know that 0 is connected to distance $n$, if none of them do, then the converse is true.  {\bf Right.} The $S$ shape of the function $p\mapsto {\bf H}_{n+1,n}(p)$. The function goes very quickly from $\varepsilon$ to $1-\varepsilon$ ($\Delta p(\varepsilon)$ is small). }\end{center}
\end{figure}

\subsection{Higher dimensions}

\paragraph{Bernoulli percolation} Again, let us start with the discussion of this simpler case. When working in higher dimensions, one can still consider crossing probabilities of boxes but one is soon facing some substantial challenges. Of course, some of the arguments of the previous section survive. For instance, one can adapt the two-dimensional proof to show that if the box $[0,n]\times[0,2n]^{d-1}$ is crossed from left to right with probability smaller than some constant $\varepsilon=\varepsilon(d)>0$, then $\theta_n(p)$ decays exponentially fast. One can also prove the differential inequality \eqref{eq:a3} without much trouble. But that is basically it.
One cannot prove that, if the probability of the box $[0,2n]\times[0,n]^{d-1}$ is crossed from left to right with probability close to one, then $\theta(p)>0$. Summarizing, we cannot (yet) exclude a regime of values of $p$ for which crossing probabilities tend to 1 but the probability that there exists an infinite connected component is zero\footnote{It is in fact the case that for $p=p_c$ and $d\ge6$, crossing probabilities tend to 1 but $\theta(p_c)=0$. What we wish to exclude is a whole range of parameters for which this happens.}.

We are therefore pushed to abandon crossing probabilities and try to work directly with $\theta_n$. Applying the BKKKL result to $\theta_n$ implies, when $p<p_c$, an inequality (basically) stating
$$\theta_n'\ge c\log n\cdot\theta_n(1-\theta_n).$$
This differential inequality is unfortunately not useful to exclude a regime where $\theta_n$ would decay polynomially fast. For this reason, we need to strengthen it. 
In order to do this, we will not rely on a concentration inequality coming from discrete Fourier analysis like in the two-dimensional case, but rather on another concentration-type inequality used in computer science. 

Informally speaking, a {\em randomized algorithm} associated with a boolean function ${\mathbf f}$ takes $\omega\in\{0,1\}^\mathbb E$ as an input, and reveals algorithmically the value of $\omega$ at different edges one by one until the value of ${\mathbf f}(\omega)$ is determined. At each step, which edge will be revealed next depends on the values of $\omega$ revealed so far. The algorithm stops as soon as the value of ${\mathbf f}$ is the same no matter the values of $\omega$ on the remaining coordinates. %Then, the question is often to determine how many bits of information must be revealed before the algorithm stops (this quantity is sometimes referred to as the computational complexity of the boolean function).

The OSSS inequality, originally introduced by O'Donnell, Saks, Schramm and Servedio in \cite{OSSS} as a step toward a conjecture of Yao \cite{yao1977probabilistic}, relates the variance of a boolean function to the influence and the computational complexity of a randomized algorithm for this function.  More precisely,
consider $p\in[0,1]$ and $n\in\mathbb N$. Fix an increasing boolean function ${\mathbf f}:\{0,1\}^\mathbb E\longrightarrow \{0,1\}$ and an algorithm ${\mathbf T}$ for ${\bf f}$. We have
\begin{equation}
    \label{eq:OSSS}\tag{OSSS}
 \mathrm{Var}_p({\mathbf f})~\le~   2 \sum_{e\in\mathbb E}  \delta_e({\mathbf T}) \, {\rm Cov}_p[{\bf f},\omega(e)],
  \end{equation}
  where 
$
\delta_e({\mathbf T})
$
is the probability that the edge $e$ is revealed by the algorithm before it stops. One will note the similarity with \eqref{eq:a1}, where the term $\delta_e(\mathbf T)$ replaces $-1$ divided by the logarithm of the covariance.

The interest of \eqref{eq:OSSS} comes from the fact that, if there exists a randomized algorithm for ${\bf f}=\mathbbm1[0\text{ connected to distance $n$}]$ for which each edge has small probability of being revealed, then the inequality implies that the derivative of $\mathbb E_p[\mathbf f]$ is much larger than the variance $\theta_n(1-\theta_n)$ of $\mathbf f$. Of course, there are several possible choices for the algorithm. Using the one described in Fig.~\ref{fig:algo}, one deduces that the probability of being revealed is bounded by $cS_n/n$  uniformly for every edge, where $S_n:=\sum_{k=0}^{n-1}\theta_k$. We therefore deduce an inequality of the form
\begin{equation}\label{eq:a4}\theta_n'\ge c'\tfrac{n}{S_n}\theta_n(1-\theta_n).
\end{equation}
Note that the quantity $n/S_n(p)$ is large when the values $\theta_k(p)$ are small, which is typically the case when $p<p_c$. In particular, one can use this differential inequality to prove the sharpness of the phase transition on any transitive graph.
Equation~\eqref{eq:a4} already appeared in Menshikov's 1986 proof while Aizenman and Barsky \cite{AizBar87} and later \cite{DumTas15} invoked alternative differential inequalities. 

\paragraph{FK percolation} As mentioned in the previous paragraph, the use of differential inequalities to prove sharpness of the phase transition is not new, and even the differential inequality \eqref{eq:a4} chosen above already appeared in the literature. Nonetheless, the existing proofs of these differential inequalities all had one feature in common: they relied on a special correlation inequality for Bernoulli percolation known as the BK inequality, which is not satisfied for most dependent percolation models, so that the historical proofs did not extend easily to FK percolation, contrarily to the approach using the OSSS inequality proposed in the previous section.

Indeed, while the OSSS inequality uses independence, it does not rely on it in a substantial way. In particular, the OSSS inequality can be extended to FK percolation, very much like \cite{GraGri06} generalized \eqref{eq:a2}. This generalization enables one to show \eqref{eq:a4} for a large class of models including dependent percolation models or so-called continuum percolation models \cite{DumRaoTas17a,DumRaoTas17c}. In particular,
\begin{theorem}[\cite{DumRaoTas17}]
Fix $q\ge1$ and $d\ge2$. Consider FK percolation on $\mathbb Z^d$ with cluster-weight $q\ge1$. For any $p<p_c$, there exists $c=c(p,q)>0$ such that for every $n\ge1$,
$$\theta_n(p,q)\le \exp(-c n).$$
\end{theorem}

Exactly as for Bernoulli percolation, one can prove many things using this theorem. Of special interest are the consequences for the Potts model (and its special case the Ising model): the exponential decay of $\theta_n(p,q)$ implies the exponential decay of correlations in the disordered phase. 
\medbreak
{\em The story of the proof of exponential decay of $\theta_n(p,q)$  is typical of percolation. Some proofs first appeared for Bernoulli percolation. These proofs were then made more robust using some external tools, here discrete analysis (the BKKKL concentration inequality or the OSSS inequality), and finally extended to more general percolation models. The next section provides another example of such a succession of events.}

\section{Computation of critical points in two dimensions}
It is often convenient to have an explicit formula for the critical point of a model. In general, one cannot really hope for such a formula but in some cases, one is saved by specific properties of the model, which can be of (at least) two kinds: {\em self-duality} or {\em exact integrability}.
 
\subsection{Computation of the critical point using self-duality}

\paragraph{Bernoulli percolation} One can easily guess why the critical point of Bernoulli percolation on $\mathbb Z^2$ should be equal to 1/2. Indeed, every configuration $\omega$ is naturally associated with a dual configuration $\omega^*$ defined on the dual lattice $(\mathbb Z^2)^*=(\tfrac12,\tfrac12)+\mathbb Z^2$ of $\mathbb Z^2$: for every edge $e$, set 
$$\omega^*(e^*)\stackrel{\rm def}=1-\omega(e),$$ where $e^*$ is the unique edge of the dual lattice crossing the edge $e$ in its middle. In words, a dual edge is open if the corresponding edge of the primal lattice is closed, and vice versa. If $\omega$ is sampled according to Bernoulli percolation of parameter $p$, then $\omega^*$ is sampled according to a Bernoulli percolation on $(\mathbb Z^2)^*$ of parameter $p^*:=1-p$. The value $1/2$ therefore emerges as the self-dual value for which $p=p^*$. 

It is not a priori clear why the self-dual value should be the critical one, but armed with the theorems of the previous section, we can turn this observation into a rigorous proof. Indeed, one may check (see Fig.~\ref{fig:6}) that for every $n\ge1$,
$${\bf H}_{n+1,n}(\tfrac12)=\tfrac12.$$
Yet, an outcome of Section~\ref{sec:3} is that crossing probabilities are tending to 0 when $p<p_c$ and to 1 when $p>p_c$. As a consequence, the only possible value for $p_c$ is $1/2$. 

\paragraph{FK percolation}The duality relation generalizes to cluster-weights $q\ne1$: if $\omega$ is sampled according to a FK percolation measure with parameters $p$ and $q$, then $\omega^*$ is sampled according to a FK percolation measure with parameters $p^*$ and $q^*$ satisfying 
$$\frac{pp^*}{(1-p)(1-p^*)}=q\qquad\text{and}\qquad q^*=q.$$
The proof of this fact involves Euler's relation for planar graphs. Let us remark that readers trying to obtain such a statement as an exercise will encounter a small difficulty due to boundary effects on $\mathbb G$; we refer to \cite{Dum17} for details how to handle such {\em boundary conditions}. 
The formulas above imply that for every $q\ne0$, there exists a unique point $p_{\rm sd}(q)$ such that 
$$p_{\rm sd}(q)=p_{\rm sd}(q)^*=\frac{\sqrt q}{1+\sqrt q}.$$
Exactly as in the case of Bernoulli percolation, one may deduce from self-duality some estimates on crossing probabilities at $p=p_{\rm sd}(q)$, which imply in the very same way the following theorem.
\begin{theorem}[\cite{BefDum12}]
The critical point of FK percolation on the square lattice with cluster-weight $q\ge1$ is equal to the self-dual point $\sqrt q/(1+\sqrt q)$.
\end{theorem}

\subsection{Computation via parafermionic observables}

Sometimes, no obvious self-duality relation helps us identify the critical point, but one can be saved by a second strategy. In order to illustrate it, consider the loop $O(n)$ model with parameters $x>0$ and $n\in[0,2]$ and its associated FK representation described in Section~\ref{sec:2}. Rather than referring to duality (in this case, none is available as for today), the idea consists in introducing a function that satisfies some specific integrability/local relations at a given value of the parameter.

\begin{figure}
\begin{center}\includegraphics[width=0.25\textwidth]{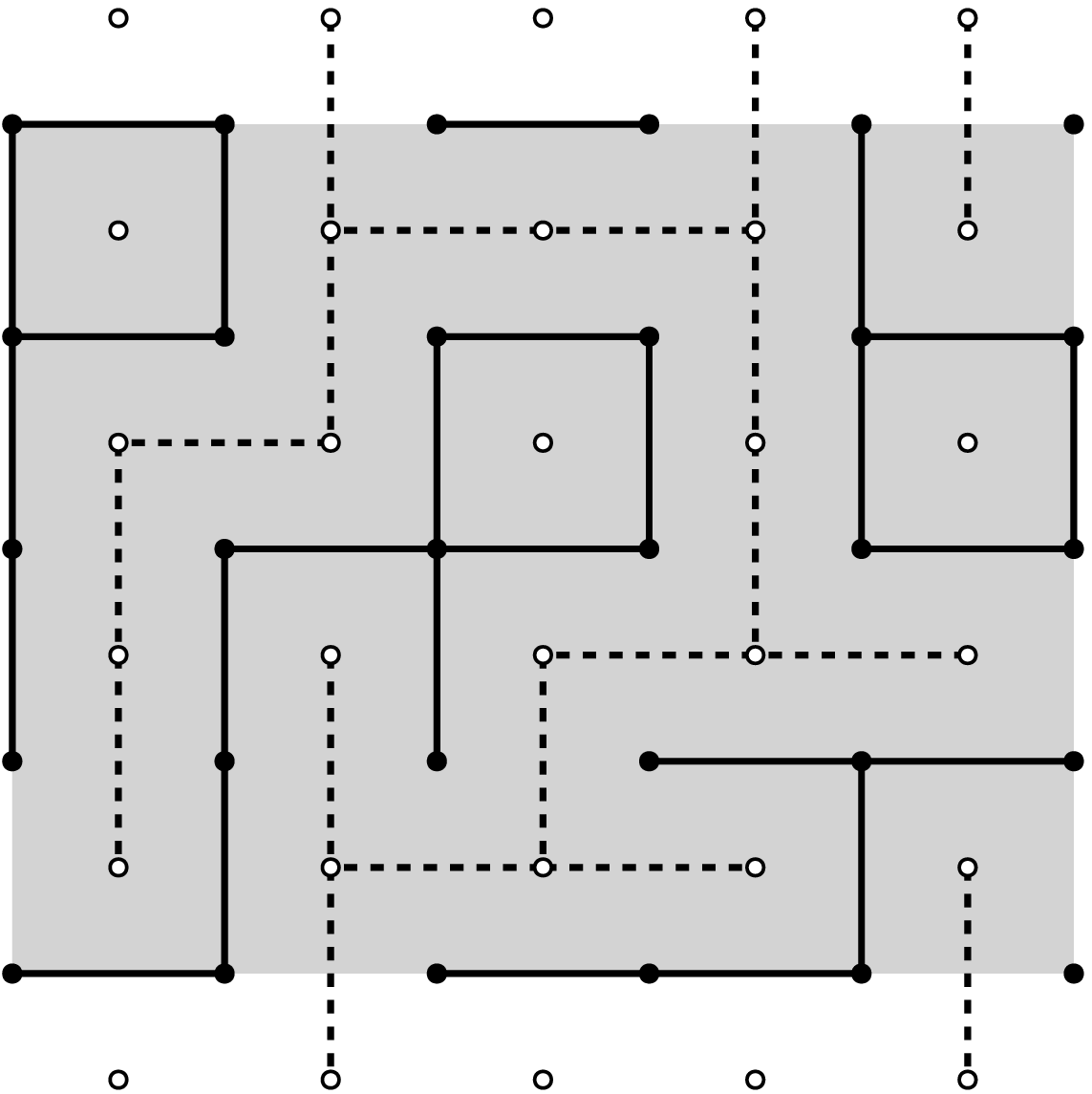}\qquad\qquad\qquad\includegraphics[width=0.45\textwidth]{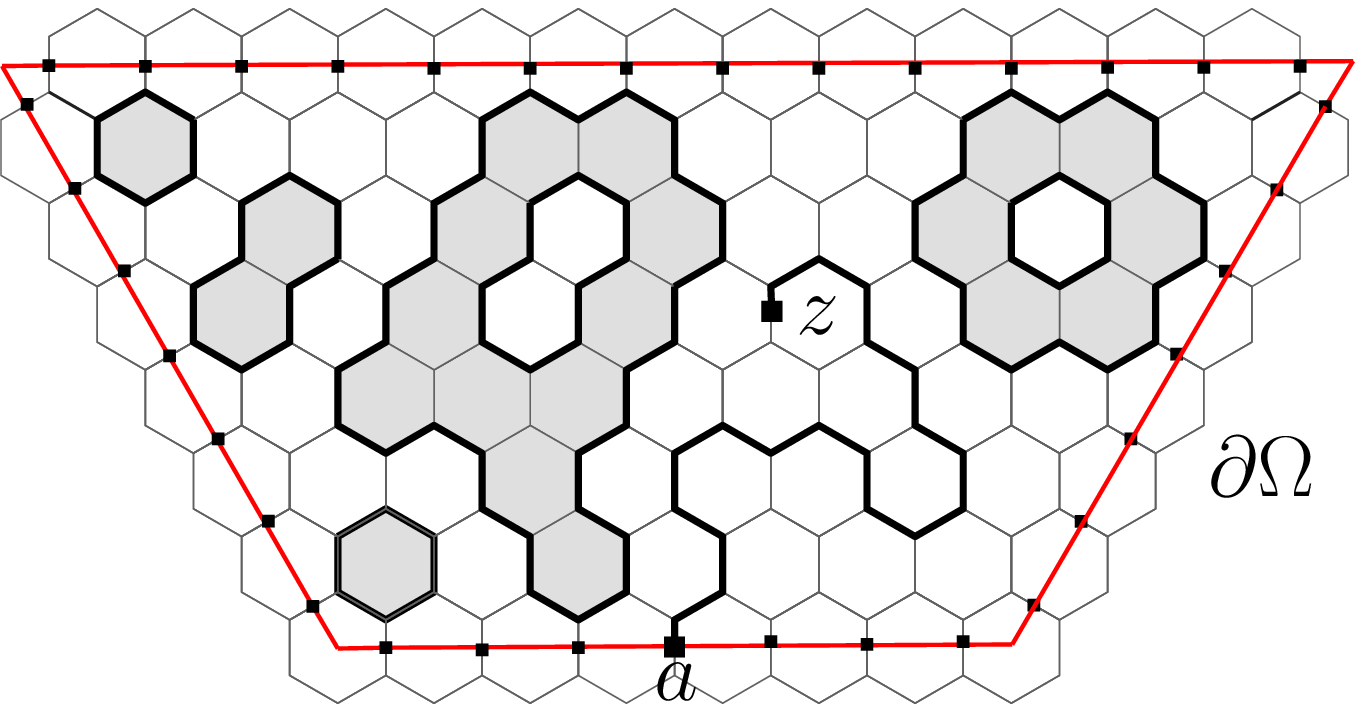}
\caption{\label{fig:6}{\bf Left.} If $[0,n+1]\times[0,n]$ is not crossed from left to right, then the boundary of the connected components touching the right side is a dual path going from top to bottom. {\bf Right.} The domain $\Omega$ with its boundary $\partial\Omega$. The configuration $\omega$ corresponds to the interfaces between the gray and the white hexagons. The path $\gamma$ runs from $a$ to $z$ without intersecting $\omega$. }\end{center}
\end{figure}

Take a self-avoiding polygon on the dual (triangular) lattice of the hexagonal lattice; see Fig.~\ref{fig:6}. By definition, this polygon divides the hexagonal lattice into two connected components, a bounded one and an unbounded one. Call the bounded one $\Omega$ and, by analogy with the continuum, denote the self-avoiding polygon by $\partial\Omega$.

Define the  \emph{parafermionic observable} introduced in \cite{Smi11}, as follows (see Fig.~\ref{fig:6}): for a mid-edge $z$ in $\Omega$ and a mid-edge $a$ in $\partial\Omega$, set
\begin{equation*}
	F(z)=F(\Omega,a,z,n,x,\sigma)\stackrel{\rm def}=\sum_{\substack{\omega,\gamma\subset \Omega\\ \omega\cap\gamma=\emptyset}}{\rm e}^{-{\rm i}\sigma \mathbf W_\gamma(a,z)} x^{|\gamma|+|\omega|}n^{\ell(\omega)}
\end{equation*}
(recall that $\ell(\omega)$ is the number of loops in $\omega$), where the sum is over pairs $(\omega,\gamma)$ with $\omega$ a loop configuration, and $\gamma$ a self-avoiding path from $a$ to $z$. The quantity $\mathbf W_\gamma(a,z)$, called the {\em winding term}, is equal to $\tfrac\pi3$ times the number of left turns minus the number of right turns  made by the walk $\gamma$ when going from $a$ to $z$. It corresponds to the total rotation of the oriented path $\gamma$.

The interest of $F$ lies in a special property satisfied when the parameters of the model are tuned properly. More precisely, if $\sigma=\sigma(n)$ is well chosen (the formula is explicit but irrelevant here) and
$$x=x_c(n)\stackrel{\rm def}=\frac{1}{\sqrt{2+\sqrt{2-n}}},$$ the function satisfies that for any self-avoiding polygon $\mathbf C=(c_0,c_1,\dots,c_k=c_0)$ on $\mathbb T$ that remains within the bounded region delimited by $\partial\Omega$,
\begin{equation}\label{eq:j}\oint_{\mathbf C}F(z)dz\stackrel{\rm def}=\sum_{i=1}^k(c_i-c_{i-1})F(\tfrac{c_{i-1}+c_i}2)=0.\end{equation}
Above, the quantity $c_i$ is considered as a complex number, in such a way that the previous definition matches the intuitive notion of contour integral for functions defined on the middle of the edges of the hexagonal lattice.

In words, \eqref{eq:j} means that for special values of $x$ and $\sigma$, discrete contour integrals of the parafermionic observable vanish. 
	In the light of Morera's theorem, this property is a glimpse of conformal invariance of the model in the sense that the observable satisfies a weak notion of discrete holomorphicity. This singles out $x_c(n)$ as a very peculiar value of the parameter $x$.
	
	The existence of such a discrete holomorphic observable at $x_c(n)$ did not really come as a surprise. In the case of the loop $O(n)$ model with loop-weight $n\in[0,2]$, the physicist Nienhuis \cite{Nie82} predicted that $x_c(n)$ is a critical value for the loop model, in the following sense: 
	\begin{itemize}
	\item when $x<x_c(n)$, loops are typically small: the probability that the loop of the origin is of length $k$ decays exponentially fast in $k$. 
	  \item when $x\ge x_c(n)$, loops are typically large: the probability that the loop of the origin is of length $k$ decays slower than polynomially fast in $k$.
  \end{itemize}
Physically, this criticality of $x_c(n)$ has an important consequence. As (briefly) mentioned in Section~\ref{sec:11}, Bernoulli percolation and more generally two-dimensional models at criticality are predicted to be conformally invariant. This prediction has a concrete implication on critical models: some observables\footnote{Roughly speaking, observables are averages of random quantities defined locally in terms of the system.} should converge in the scaling limit to conformally invariant/covariant objects. In the continuum, typical examples of such objects are provided by harmonic and holomorphic solutions to Boundary Value Problems; it is thus natural to expect that some observables of the loop $O(n)$ model at criticality are discrete holomorphic. 
	
	\begin{remark}Other than being interesting in themselves, discrete holomorphic
functions have found several applications in geometry, analysis, combinatorics, and probability. The use of discrete holomorphicity in statistical physics
appeared first in the case of dimers \cite{Ken00} and has since then been extended to several statistical
physics models, including the Ising model \cite{CheSmi11,CheSmi12} (see also \cite{DumSmi12a} for references). \end{remark}

The definition of discrete holomorphicity usually imposes stronger conditions on the function than just zero contour integrals (for instance, one often asks for suitable discretizations of the Cauchy-Riemann equations).
In particular, in our case the zero contour integrals do not uniquely determine $F$. Indeed, there is one variable $F(z)$ by edge, but the number of independent linear relations provided by the zero contour integrals is equal to the number of vertices\footnote{Indeed, it is sufficient to know that discrete contour integrals vanish for a basis of the $\mathbb Z$-module of cycles (which are exactly the contours) on the triangular lattice staying in $\Omega$ to obtain all the relations in \eqref{eq:j}. A natural choice for such a basis is provided by the triangular cycles around each face of the triangular lattice inside $\partial\Omega$, hence it has exactly as many elements as vertices of the hexagonal lattice in $\Omega$.}. In conclusion, there are much fewer relations than unknown and it is completely unclear whether one can extract any relevant information from \eqref{eq:j}.

Anyway, one can still try to harvest this ``partial'' information to identify rigorously the critical value of the loop $O(n)$ model. For $n=0$ (in this case there is no loop configuration and just one self-avoiding path), the parafermionic observable was used to show that the {\em connective constant} of the hexagonal lattice \cite{DumSmi12}, i.e.
$$\mu_c\stackrel{\rm def}=\lim_{n\rightarrow\infty} \#\{\text{self-avoiding walks of length $n$ starting at the origin}\}^{1/n}$$
 is equal to $\sqrt{2+\sqrt2}$. For $n\in[1,2]$, the same observable was used to show that at $x=x_c(n)$, the probability of having a loop of length $k$ decays slower than polynomially fast \cite{DumGlaPelSpi17}, thus proving part of Nienhuis prediction (this work has also applications for the corresponding site percolation model described in Section~\ref{sec:2}). 

Let us conclude this section by mentioning that the parafermionic observable defined for the loop $O(n)$ model can also be defined for a wide variety of models of planar statistical physics; see e.g.~\cite{IkhCar09,IkhWesWhe13,RajCar07}. This leaves hope that many more models from planar statistical physics can be studied using discrete holomorphic functions. For the FK percolation of parameter $q\ge4$, a parafermionic observable was used to show that $p_c(q)=\sqrt q/(+\sqrt q)$, thus providing an alternative proof to \cite{BefDum12} (this proof was in fact obtained prior to the proof of \cite{BefDum12}). Recently, the argument was generalized to the case $q\in[1,4]$ in \cite{MukSmi17}.
\medbreak
{\em In the last two sections, we explained how the study of crossing probabilities can be combined with duality or parafermionic observables to identify the critical value of some percolation models. In the next one, we go further and discuss how the same tools can be used to decide whether the phase transition is continuous or discontinuous (see definition below).
}
\section{The critical behavior of planar dependent percolation models}

\subsection{Renormalization of crossing probabilities}

For Bernoulli percolation, we mentioned that crossing probabilities remain bounded away from 0 and 1 uniformly in the size of the rectangle (provided the aspect ratio stays bounded away from 0 or 1). For more complicated percolation models, the question is more delicate, in particular due to the long-range dependencies. For instance, it may be that crossing probabilities tend to zero when conditioning on edges outside the box to be closed, and to 1 if these edges are conditioned to be open.
To circumvent this problem, we introduce a new property.

Consider a percolation measure $\mathbb P$ (one can think of a FK percolation measure for instance). Let $\Lambda_n$ be the box of size $n$ around the origin. We say that $\mathbb P$ satisfies the {\em (polynomial) mixing property} if
\medbreak\noindent{\bf (Mix)} {\em there exist $c,C\in(0,\infty)$ such that for every $N\ge 2n$ and every events $A$ and $B$ depending on edges in $\Lambda_n$ and outside $\Lambda_N$ respectively, we have that }
$$|\mathbb P[A\cap B]-\mathbb P[A]\mathbb P[B]|\le C(\tfrac nN)^c\cdot\mathbb P[A]\mathbb P[B].$$
This property has many implications for the study of the percolation model, mostly because it enables one to decorrelate events happening in different parts of the space. 

It is a priori unclear how one may prove the mixing property in general. Nonetheless, for critical FK percolation, it can be shown that (wMix) is equivalent to the {\em strong box crossing property}: uniformly on the states of edges outside of $\Lambda_{2n}$, crossing a rectangle of aspect ratio $\rho$ included in $\Lambda_n$ remains bounded away from 0 and 1 uniformly in $n$. Note that the difference with the previous sections comes from the fact that we consider crossing probabilities conditioned on the state of edges at distance $n$ of the rectangle (of course, when considering Bernoulli percolation, this does not change anything, but this is not the case anymore when $q>1$).
 
The mixing property is not always satisfied at criticality. Nevertheless, in \cite{DumSidTas13}
the following dichotomy result was obtained.
\begin{theorem}[The continuous/discontinuous dichotomy]
  For any $q\ge1$, 
  \begin{itemize}
  \item either (wMix) is satisfied. In such case:
  \begin{itemize}[noitemsep,nolistsep]
  \item $\theta(p,q)$ tends to 0 as $p\searrow p_c$;
  \item There exists $c>0$ such that $cn^{-1}\le \theta_n(p_c,q)\le n^{-c}.$
  \item Crossing probabilities of a rectangle of size roughly $n$ remain bounded away from 0 and 1 uniformly in the state of edges at distance $n$ of the rectangle;
%  \item There exists a unique infinite-volume FK percolation measure;
  \item The rate of exponential decay of $\theta_n(p,q)$ goes to 0 as $p\nearrow p_c$.
  \end{itemize}
  \item or (wMix) is not satisfied and in such case:
  \begin{itemize}[noitemsep,nolistsep]
  \item $\theta(p,q)$ does not tend to 0 as $p\searrow p_c$;
  \item There exists $c>0$ such that for every $n\ge1$,
  $\theta_n(p_c,q)\le \exp(-cn)$;
  \item Crossing probabilities of a rectangle of size roughly $n$ tend to 0 (resp.~1) when conditioned on the state of edges at distance $n$ of the rectangle to be closed (resp.~open);
%  \item There are distinct infinite-volume FK percolation measures;
  \item The rate of exponential decay of $\theta_n(p,q)$ does not go to 0 as $p\nearrow p_c$.
  \end{itemize}
  \end{itemize}
\end{theorem}

In the first case, we say that the phase transition is {\em continuous} in reference to the fact that $p\mapsto \theta(p,q)$ is continuous at $p_c$. In the second case, we say that the phase transition is {\em discontinuous}. Interestingly, this result also shows that a number of (potentially different) definitions of continuous/discontinuous phase transitions sometimes used  in physics are in fact the same one in the case of FK percolation.

The proof of the dichotomy is based on a renormalization scheme for crossing probabilities when conditioned on edges outside a box to be closed. Explaining the strategy would lead us too far, and we refer to \cite{DumRaoTas17b,DumSidTas13} for more details. Let us simply add that the proof is not specific to FK percolation and has been extended to other percolation models (see for instance \cite{DumGlaPelSpi17,DumRaoTas17b}), so one should not think of this result as an isolated property of FK percolation, but rather as a general feature of two-dimensional dependent percolation models.

% for the case of the site percolation model obtained from the loop $O(n)$ model with loop-weight $n\ge1$ and edge-weight $x$ satisfying $x^2n\le 1$. 
%\begin{remark}The dichotomy result combines with the loop $O(n)$ result from the last section, it does imply that the model has large loops for $x=x_c(n)$ and $n\in[1,2]$, as predicted by Nienhuis.\end{remark}
\begin{figure}
\begin{center}\includegraphics[width=0.23\textwidth]{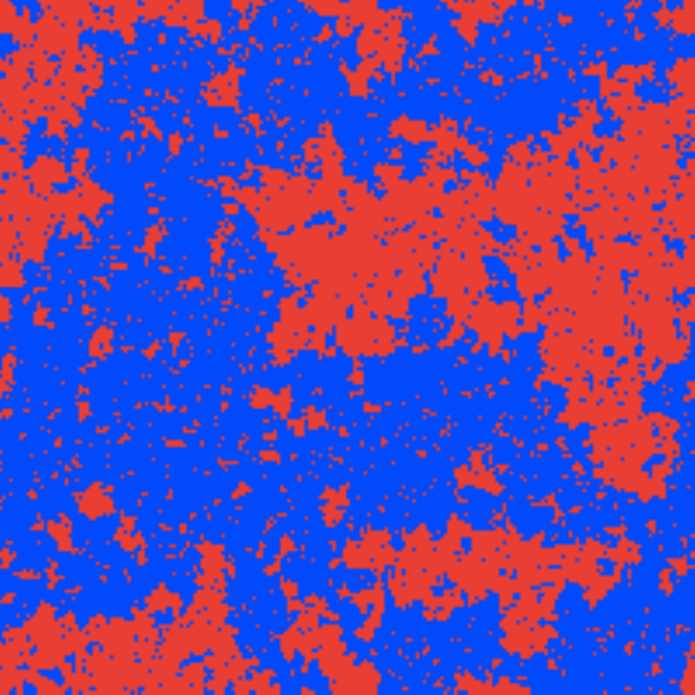}\quad\includegraphics[width=0.23\textwidth]{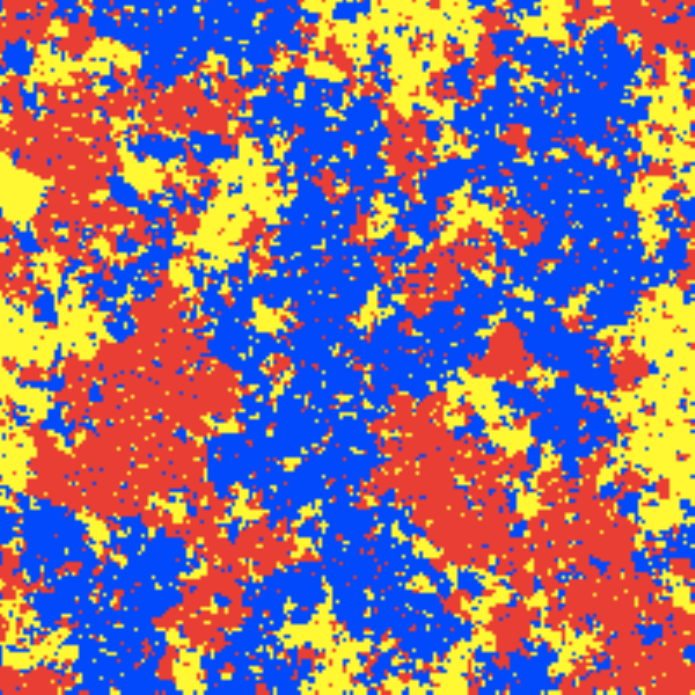}\quad\includegraphics[width=0.23\textwidth]{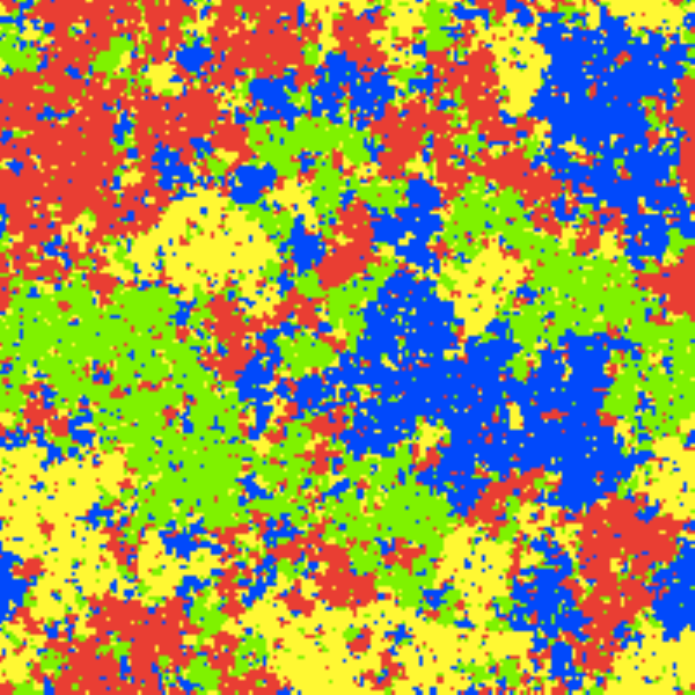}
\quad\includegraphics[width=0.23\textwidth]{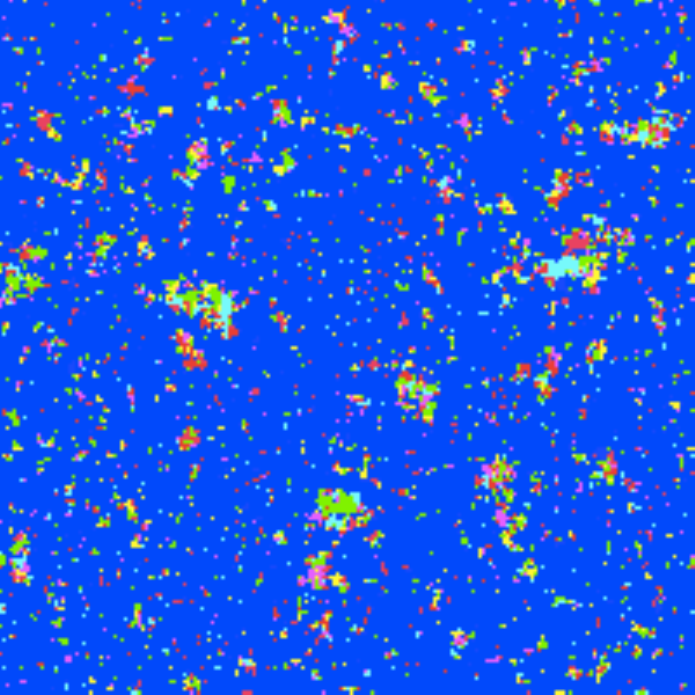}
\caption{\label{fig:10}Simulations, due to Vincent Beffara, of the critical planar Potts model with $q$ equal to $2$, $3$, $4$, and $9$ respectively. The behavior for $q\le 4$ is clearly different from the behavior for $q=9$. In the first three pictures, each color seems to play the same role, while in the last three, one color wins over the other ones. This is coherent since the phase transition of the associated FK percolation is continuous for $q\le 4$ and discontinuous for $q\ge 5$.}\end{center}
\end{figure}
\subsection{Deciding the dichotomy} 

As mentioned above, critical planar dependent percolation models can exhibit two different types of critical behavior: continuous or discontinuous. In order to decide which one of the two it is, one needs to work a little bit harder. Let us (briefly) describe two possible strategies. We restrict to the case of FK percolation, for which Baxter \cite{Bax73} conjectured that for $q\le q_c(2)=4$, the phase transition is continuous, and for $q>q_c(2)$, the phase transition is discontinuous; see Fig.~\ref{fig:10}.

To prove that the phase transition is discontinuous for $q>4$, we used  a method going back to early works on the six-vertex model \cite{DumGanHar16}. 
The six-vertex model was initially proposed by Pauling in 1931 in order to study the thermodynamic properties of ice.
While we are mainly interested in its connection to the previously discussed models, 
the six-vertex model is a major object of study on its own right; we refer to \cite{Resh10} and Chapter~8 of \cite{Bax89} (and references therein) for the definition and a bibliography on the subject. 

The utility of the six-vertex model stems from its solvability using the transfer-matrix formalism. More precisely, the partition function of a six-vertex model on a torus of size $N$ times $M$ may be expressed as the trace of the $M$-th power of a matrix $V$ (depending on $N$) called the {\em transfer matrix}. This property can be used to   rigorously compute the Perron-Frobenius eigenvalues of the diagonal blocks of the transfer matrix,
whose ratios are then related to the rate of exponential decay $\tau(q)$ of $\theta_n(p_c,q)$. The explicit formula obtained for $\tau(q)$ is then proved to be strictly positive for $q>4$. We should mention that this strategy is extensively used in the physics literature, in particular in the fundamental works of Baxter (again, we refer to \cite{Bax89}).

In order to prove that the phase transition is continuous for $q\le 4$, one may use the same strategy and try to prove that $\tau(q)$ is equal to 0. Nevertheless, this does not seem so simple to do rigorously, so that we prefer an alternative approach. The fact that discrete contour integrals of the parafermionic observable vanish can be used for more than just identifying the critical point. For $q\in[1,4]$, it in fact implies lower bounds on $\theta_n(p_c,q)$. These lower bounds decay at most polynomially  fast, thus guaranteeing that the phase transition is continuous thanks to the dichotomy result of the previous section. This strategy was implemented in \cite{DumSidTas13} to complete the proof of Baxter's prediction regarding the continuity/discontinuity of the phase transition for the planar FK percolation with $q\ge1$.

Let us conclude this short review by mentioning that for the special value of $q=2$, the parafermionic observable satisfies stronger constraints. This was used to show that, for this value of $q$, the observable is conformally covariant in the scaling limit  \cite{Smi10} (this paper by Smirnov had a resounding impact on our understanding of FK percolation with $q=2$), that the strong RSW property is satisfied \cite{DumHonNol11}, and that interfaces converge \cite{BenHon16,CheDumHon14}. This could be the object of a review on its own, especially regarding the conjectures generalizing these results to other values of $q$, but we reached the end of our allowed space. We refer to \cite{DumSmi12a} for more details and relevant references.

\begin{remark}
The strategy described above is very two-dimensional in nature since it relies on planarity in several occasions (crossing probabilities for the dichotomy result, parafermionic observables or transfer matrix formalism for deciding between continuity or discontinuity). In higher dimensions, the situation is more challenging. We have seen that even for Bernoulli percolation, continuity of $\theta(p)$ had not yet been proved for dimensions $3\le d\le 10$. Let us mention that for FK percolation, several results are nevertheless known. One can prove the continuity of $\theta(p,2)$  \cite{AizDumSid15} using properties specific to the Ising model (which is associated with the FK percolation with cluster-weight $q=2$ via the Edwards-Sokal coupling). Using the mean-field approximation and Reflection-Positivity, one may also show that the phase transition is discontinuous  if $d$ is fixed and $q\ge q_c(d)\gg1$ \cite{KotShl82}, 
or if $q\ge3$ is fixed and $d\ge d_c(q)\gg1$ \cite{BisCha03}. The conjecture that $q_c(d)$ is equal to 2 for any $d\ge3$ remains widely open and represents a beautiful challenge for future mathematical physicists. \end{remark}

\footnotesize
\bibliographystyle{siam}

\begin{thebibliography}{10}

\bibitem{AizBar87}
{\sc M.~Aizenman and D.~J. Barsky}, {\em Sharpness of the phase transition in
  percolation models}, Comm. Math. Phys., 108 (1987), pp.~489--526.

\bibitem{AizDumSid15}
{\sc M.~Aizenman, H.~{Duminil-Copin}, and V.~Sidoravicius}, {\em Random
  {C}urrents and {C}ontinuity of {I}sing {M}odel's {S}pontaneous
  {M}agnetization}, Communications in Mathematical Physics, 334 (2015),
  pp.~719--742.

\bibitem{AizKesNew87}
{\sc M.~Aizenman, H.~Kesten, and C.~M. Newman}, {\em Uniqueness of the infinite
  cluster and continuity of connectivity functions for short and long range
  percolation}, Comm. Math. Phys., 111 (1987), pp.~505--531.

\bibitem{BarGriNew91}
{\sc D.~J. Barsky, G.~R. Grimmett, and C.~M. Newman}, {\em Percolation in
  half-spaces: equality of critical densities and continuity of the percolation
  probability}, Probab. Theory Related Fields, 90 (1991), pp.~111--148.

\bibitem{Bax73}
{\sc R.~J. Baxter}, {\em Potts model at the critical temperature}, Journal of
  Physics C: Solid State Physics, 6 (1973), p.~L445.

\bibitem{Bax89}
{\sc R.~J. Baxter}, {\em Exactly solved models in statistical mechanics},
  Academic Press Inc. [Harcourt Brace Jovanovich Publishers], London, 1989.
\newblock Reprint of the 1982 original.

\bibitem{Bek75}
{\sc W.~Beckner}, {\em Inequalities in Fourier analysis}, Ann. of Math, 102
  (1975), pp.~159--182.

\bibitem{BefDum12}
{\sc V.~Beffara and H.~{Duminil-Copin}}, {\em The self-dual point of the
  two-dimensional random-cluster model is critical for {$q\geq 1$}}, Probab.
  Theory Related Fields, 153 (2012), pp.~511--542.

\bibitem{BefDum13}
{\sc V.~Beffara and H.~{Duminil-Copin}}, {\em Lectures on planar
  percolation with a glimpse of {S}chramm {L}oewner {E}volution}, Probability
  Surveys, 10 (2013), pp.~1--50.

\bibitem{BelPolZam84}
{\sc A.~A. Belavin, A.~M. Polyakov, and A.~B. Zamolodchikov}, {\em Infinite
  conformal symmetry of critical fluctuations in two dimensions}, J. Statist.
  Phys., 34 (1984), pp.~763--774.

\bibitem{BenHon16}
{\sc S.~Benoist and C.~Hongler},
{\em The scaling limit of critical ising interfaces is {CLE}(3)},
arXiv:1604.06975.

\bibitem{BenKalSch99}
{\sc I.~Benjamini, G.~Kalai, and O.~Schramm}, {\em Noise sensitivity of
  {B}oolean functions and applications to percolation}, Inst. Hautes {\'E}tudes
  Sci. Publ. Math.,  (1999), pp.~5--43.

\bibitem{BenSch96}
{\sc I.~Benjamini and O.~Schramm}, {\em Percolation beyond {$\bf Z^d$}, many
  questions and a few answers}, Electron. Comm. Probab., 1 (1996), pp.~no.\ 8,
  71--82 (electronic).

\bibitem{BisCha03}
{\sc M.~Biskup and L.~Chayes}, {\em Rigorous analysis of discontinuous phase
  transitions via mean-field bounds}, Comm. Math. Phys,  (2003), pp.~53--93.

\bibitem{BolRio06}
{\sc B.~Bollob{{\'a}}s and O.~Riordan}, {\em The critical probability for
  random {V}oronoi percolation in the plane is 1/2}, Probab. Theory Related
  Fields, 136 (2006), pp.~417--468.

\bibitem{BolRio06c}
{\sc B.~Bollob{{\'a}}s and O.~Riordan}, {\em A short proof of the
  {H}arris-{K}esten theorem}, Bull. London Math. Soc., 38 (2006), pp.~470--484.

\bibitem{BolRio10}
{\sc B.~Bollob{{\'a}}s and O.~Riordan}, {\em Percolation on
  self-dual polygon configurations}, in An irregular mind, vol.~21 of Bolyai
  Soc. Math. Stud., J{\'a}nos Bolyai Math. Soc., Budapest, 2010, pp.~131--217.

\bibitem{Bon70}
{\sc A.~Bonami}, {\em Etude des coefficients de Fourier des fonctions de
  lp(g)}, Ann. Inst. Fourier, 20 (1970), pp.~335--402.

\bibitem{BouKahKal92}
{\sc J.~Bourgain, J.~Kahn, G.~Kalai, Y.~Katznelson, and N.~Linial}, {\em The
  influence of variables in product spaces}, Israel J. Math., 77 (1992),
  pp.~55--64.

\bibitem{BroHam57}
{\sc S.~R. Broadbent and J.~M. Hammersley}, {\em Percolation processes. {I}.
  {C}rystals and mazes}, Proc. Cambridge Philos. Soc., 53 (1957), pp.~629--641.

\bibitem{BurKea89}
{\sc R.~M. Burton and M.~Keane}, {\em Density and uniqueness in percolation},
  Comm. Math. Phys., 121 (1989), pp.~501--505.

\bibitem{CamNew06}
{\sc F.~Camia and C.~Newman}, {\em Two-dimensional critical percolation: the full scaling limit},
  Comm. Math. Phys., 268(1) (2006), pp.~1--38.

\bibitem{Car92}
{\sc J.~L. Cardy}, {\em Critical percolation in finite geometries}, J. Phys. A,
  25 (1992), pp.~L201--L206.

\bibitem{ChaMac97}
{\sc L.~Chayes and J.~Machta}, {\em Graphical representations and cluster
  algorithms I. Discrete spin systems}, Phys. A, 239 (1997), pp.~542--601.

\bibitem{CheDumHon14}
{\sc D.~Chelkak, H.~{Duminil-Copin}, C.~Hongler, A.~Kemppainen, and S.~Smirnov},
{\em Convergence of {I}sing interfaces to {S}chramm's {SLE} curves,}
C. R. Acad. Sci. Paris Math. , 352(2) (2014), pp.~157--161.

\bibitem{CheSmi11}
{\sc D.~Chelkak and S.~Smirnov}, {\em Discrete complex analysis on isoradial
  graphs}, Adv. Math., 228 (2011), pp.~1590--1630.

\bibitem{CheSmi12}
{\sc D.~Chelkak and S.~Smirnov}, {\em Universality in the
  2{D} {I}sing model and conformal invariance of fermionic observables},
  Invent. Math., 189 (2012), pp.~515--580.

\bibitem{DomMukNie81}
{\sc E.~Domany, D.~Mukamel, B.~Nienhuis, and A.~Schwimmer}, {\em Duality
  relations and equivalences for models with {O}(n) and cubic symmetry},
  Nuclear Physics B, 190 (1981), pp.~279--287.

\bibitem{Dum17}
{\sc H.~{Duminil-Copin}}, {\em Lectures on the Ising and Potts models on the hypercubic lattice}, arXiv:1707.00520, 2017.

\bibitem{DumGanHar16}
{\sc H.~Duminil-Copin, M.~Gagnebin, M.~Harel, I.~Manolescu, and V.~Tassion},
  {\em Discontinuity of the phase transition for the planar random-cluster and
  {P}otts models with $ q> 4$}, arXiv:1611.09877.

\bibitem{DumGlaPelSpi17}
{\sc H.~{Duminil-Copin}, A.~Glazman, R.~Peled, and Y.~Spinka}, {\em Macroscopic
  loops in the loop $O(n)$ model at {N}ienhuis' critical point.}
\newblock arXiv:1707.09335.

\bibitem{DumHonNol11}
{\sc H.~{Duminil-Copin}, C.~Hongler, and P.~Nolin}, {\em Connection
  probabilities and {RSW}-type bounds for the two-dimensional {FK} {I}sing
  model}, Comm. Pure Appl. Math., 64 (2011), pp.~1165--1198.

\bibitem{DumMan16}
{\sc H.~Duminil-Copin and I.~Manolescu}, {\em The phase transitions of the
  planar random-cluster and {P}otts models with $q\ge1$ are sharp}, Probability
  Theory and Related Fields, 164 (2016), pp.~865--892.

\bibitem{DumRaoTas17b}
{\sc H.~{Duminil-Copin}, A.~Raoufi, and V.~Tassion}, {\em Renormalization of
  crossings in planar dependent percolation models}.
\newblock preprint.

\bibitem{DumRaoTas16}
{\sc H.~{Duminil-Copin}, A.~Raoufi, and V.~Tassion}, {\em A new computation of
  the critical point for the planar random-cluster model with $q\ge1$}.
\newblock arXiv:1604.03702.

\bibitem{DumRaoTas17a}
{\sc H.~{Duminil-Copin}, A.~Raoufi, and V.~Tassion}, {\em Exponential decay of
  connection probabilities for subcritical {V}oronoi percolation in
  $\mathbb{R}^d$}.
\newblock arXiv:1705.07978.

\bibitem{DumRaoTas17}
{\sc H.~Duminil-Copin, A.~Raoufi, and V.~Tassion}, {\em Sharp phase transition
  for the random-cluster and {P}otts models via decision trees}.
\newblock arXiv:1705.03104.

\bibitem{DumRaoTas17c}
H.~{Duminil-Copin}, A.~Raoufi, and V.~Tassion.
\newblock Subcritical phase of $d$-dimensional {P}oisson-boolean percolation
  and its vacant set.
\newblock {\em preprint}, 2017.

\bibitem{DumSidTas13}
{\sc H.~{Duminil-Copin}, V.~Sidoravicius, and V.~Tassion}, {\em Continuity of
  the phase transition for planar random-cluster and {P}otts models with $1\le
  q\le 4$}, Communications in {M}athematical {P}hysics, 349 (2017),
  pp.~47--107.

\bibitem{DumSmi12a}
{\sc H.~{Duminil-Copin} and S.~Smirnov}, {\em Conformal invariance of lattice
  models}, in Probability and statistical physics in two and more dimensions,
  vol.~15 of Clay Math. Proc., Amer. Math. Soc., Providence, RI, 2012,
  pp.~213--276.

\bibitem{DumSmi12}
{\sc H.~{Duminil-Copin} and S.~Smirnov}, {\em The connective
  constant of the honeycomb lattice equals {$\sqrt{2+\sqrt{2}}$}}, Ann. of
  Math. (2), 175 (2012), pp.~1653--1665.

\bibitem{DumTas15}
{\sc H.~{Duminil-Copin} and V.~Tassion}, {\em A new proof of the sharpness of
  the phase transition for {B}ernoulli percolation and the {I}sing model},
  Communications in {M}athematical {P}hysics, 343 (2016), pp.~725--745.

\bibitem{ErdRen59}
{\sc P.~Erd\"os and A.~Renyi}, {\em On random graphs I}, Publicationes
  Mathematicae, 6 (1959), pp.~290--297.

\bibitem{ForKas72}
{\sc C.~M. Fortuin and P.~W. Kasteleyn}, {\em On the random-cluster model. {I}.
  {I}ntroduction and relation to other models}, Physica, 57 (1972),
  pp.~536--564.

\bibitem{GabLyo09}
{\sc D.~Gaboriau and R.~Lyons}, {\em A measurable-group-theoretic solution to
  von Neumann's problem}, Invent. Math., 177 (2009), pp.~533--540.

\bibitem{GarPetSch11}
{\sc C.~Garban, G.~Pete, and O.~Schramm}, {\em The {F}ourier spectrum of
  critical percolation}, in Selected works of {O}ded {S}chramm. {V}olume 1, 2,
  Sel. Works Probab. Stat., Springer, New York, 2011, pp.~445--530.

\bibitem{GarPetSch13}
{\sc C.~Garban, G.~Pete, and O.~Schramm}, {\em Pivotal, cluster,
  and interface measures for critical planar percolation}, J. Amer. Math. Soc.,
  26 (2013), pp.~939--1024.

\bibitem{GarPetSch17}
{\sc C.~Garban, G.~Pete, and O.~Schramm}, {\em The scaling limits of
  near-critical and dynamical percolation}, Journal of European Math Society,
  (2017).

\bibitem{GraGri06}
{\sc B.~T. Graham and G.~R. Grimmett}, {\em Influence and sharp-threshold
  theorems for monotonic measures}, Ann. Probab., 34 (2006), pp.~1726--1745.
  
  \bibitem{Gri99}
{\sc G.~R. Grimmett}, Percolation, {\em Grundlehren der Mathematischen Wissenschaften [Fundamental Principles of \mbox{Mathematical} Sciences]}, vol.~321
 (1999).
 
  \bibitem{Gri06}
{\sc G.~R. Grimmett}, The random-cluster model, {\em Grundlehren der Mathematischen Wissenschaften [Fundamental Principles of \mbox{Mathematical} Sciences]}, vol.~333
 (2006).

%\bibitem{Gri06a}
%{\sc G.~R. Grimmett and D. Welsh}, {\em John Michael Hammersley}, arXiv:0610862, 2006.
%
\bibitem{HagPerSte97}
{\sc O.~H{\"a}ggstr{\"o}m, Y.~Peres, and J.~Steif}, {\em Dynamical
  percolation}, Ann. Inst. H. Poincar\'e Probab. Statist., 33 (1997),
  pp.~497--528.

\bibitem{Ham59}
{\sc J.M.~Hammersley}, {\em Bornes sup\'erieures de la probabilit\'e critique dans un processus de filtration,
Le Calcul des Probabilit\'es et ses Applications}, CNRS, Paris (1959), pp. 17--37.

\bibitem{HarSla90}
{\sc T.~Hara and G.~Slade}, {\em Mean-field critical behaviour for percolation
  in high dimensions}, Comm. Math. Phys., 128 (1990), pp.~333--391.

\bibitem{Har60}
{\sc T.~E. Harris}, {\em A lower bound for the critical probability in a
  certain percolation process}, Proc. Cambridge Philos. Soc., 56 (1960),
  pp.~13--20.



\bibitem{IkhCar09}
{\sc Y.~Ikhlef and J.~Cardy}, {\em Discretely holomorphic parafermions and
  integrable loop models}, J. Phys. A, 42 (2009), pp.~102001, 11.

\bibitem{IkhWesWhe13}
{\sc Y.~Ikhlef, R.~Weston, M.~Wheeler, and P.~Zinn-Justin}, {\em Discrete
  holomorphicity and quantized affine algebras}.
\newblock arxiv:1302.4649.

\bibitem{KahKalLin88}
{\sc J.~Kahn, G.~Kalai, and N.~Linial}, {\em The influence of variables on
  boolean functions}, in 29th Annual Symposium on Foundations of Computer
  Science, (1988), pp.~68--80.

\bibitem{Ken00}
{\sc R.~Kenyon}, {\em Conformal invariance of domino tiling}, Ann. Probab., 28
  (2000), pp.~759--795.

\bibitem{Kes80}
{\sc H.~Kesten}, {\em The critical probability of bond percolation on the
  square lattice equals {${1\over 2}$}}, Comm. Math. Phys., 74 (1980),
  pp.~41--59.
  
  \bibitem{Kes82}
{\sc H.~Kesten}, Percolation theory for mathematicians, {\em Birkh\"auser Boston}, Progress in Probability and Statistics, vol.~2 (1982).

\bibitem{KotShl82}
{\sc R.~Koteck{{\'y}} and S.~B. Shlosman}, {\em First-order phase transitions
  in large entropy lattice models}, Comm. Math. Phys., 83 (1982), pp.~493--515.

\bibitem{LanPouSai94}
{\sc R.~Langlands, P.~Pouliot, and Y.~Saint-Aubin}, {\em Conformal invariance
  in two-dimensional percolation}, Bull. Amer. Math. Soc. (N.S.), 30 (1994),
  pp.~1--61.

\bibitem{Men86}
{\sc M.~V. Menshikov}, {\em Coincidence of critical points in percolation
  problems}, Dokl. Akad. Nauk SSSR, 288 (1986), pp.~1308--1311.

\bibitem{MukSmi17}
{\sc E.~Mukoseeva and D.~Smirnova}, {\em Computation of the critical point for
  random-cluster models via the parafermionic observable}.
\newblock preprint.

\bibitem{Nie82}
{\sc B.~Nienhuis}, {\em Exact {C}ritical {P}oint and {C}ritical {E}xponents of
  $\mathrm{O}(n)$ {M}odels in {T}wo {D}imensions}, Physical Review Letters, 49
  (1982), pp.~1062--1065.

\bibitem{OSSS}
{\sc R.~O'Donnell, M.~Saks, O.~Schramm, and R.~Servedio}, {\em Every decision
  tree has an influential variable}, FOCS,  (2005).

\bibitem{Ols80}
{\sc A.~Olsanskii}, {\em On the question of the existence of an invariant mean
  on a group}, Uspekhi Mat. Nauk, 35 (1980), pp.~199--200.

\bibitem{PakSmi00}
{\sc I.~Pak and T.~Smirnova-Nagnibeda}, {\em On non-uniqueness of percolation
  on nonamenable Cayley graphs}, C. R. Acad. Sci. Paris, 330 (2000),
  pp.~495--500.

\bibitem{PfiVel97}
{\sc C.-E. Pfister and Y.~Velenik}, {\em Random-cluster representation of the
  ashkin-teller model}, J. Stat. Phys., 88 (1997), pp.~1295--1331.

\bibitem{FitHof15}
{\sc F.~R. and van~der Hofstad~R.}, {\em Mean-field behavior for
  nearest-neighbor percolation in d>10}.
\newblock arXiv:1506.07977.

\bibitem{RajCar07}
{\sc M.~A. Rajabpour and J.~Cardy}, {\em Discretely holomorphic parafermions in
  lattice {$Z_N$} models}, J. Phys. A, 40 (2007), pp.~14703--14713.

\bibitem{Resh10}
{\sc N.~{Reshetikhin}}, {\em {Lectures on the integrability of the 6-vertex
  model}}, arXiv1010.5031.

\bibitem{Rus78}
{\sc L.~Russo}, {\em A note on percolation}, Z. Wahrscheinlichkeitstheorie und
  Verw. Gebiete, 43 (1978), pp.~39--48.

\bibitem{Rus82}
{\sc L.~Russo}, {\em An Approximate Zero-One Law }, Z. Wahrscheinlichkeitstheorie und
  Verw. Gebiete, 61 (1982), ppp.~129--139.


\bibitem{Sch00}
{\sc O.~Schramm}, {\em Scaling limits of loop-erased random walks and uniform
  spanning trees}, Israel J. Math., 118 (2000), pp.~221--288.

\bibitem{SeyWel78}
{\sc P.~D. Seymour and D.~J.~A. Welsh}, {\em Percolation probabilities on the
  square lattice}, Ann. Discrete Math., 3 (1978), pp.~227--245.
\newblock Advances in graph theory (Cambridge Combinatorial Conf., Trinity
  College, Cambridge, 1977).

\bibitem{Smi01}
{\sc S.~Smirnov}, {\em Critical percolation in the plane: conformal invariance,
  {C}ardy's formula, scaling limits}, C. R. Acad. Sci. Paris S{\'e}r. I Math.,
  333 (2001), pp.~239--244.

\bibitem{Smi10}
{\sc S.~Smirnov}, {\em Conformal invariance
  in random cluster models. {I}. {H}olomorphic fermions in the {I}sing model},
  Ann. of Math. (2), 172 (2010), pp.~1435--1467.
  
  \bibitem{Smi11}
{\sc S.~Smirnov}, {\em Discrete complex analysis and probability},
  Proceedings of the {I}nternational {C}ongress of {M}athematicians. {V}olume {I} (2010), pp.~595--621.
  
  \bibitem{Smi17}
{\sc S.~Smirnov}, private communications.

\bibitem{Tal94}
{\sc M.~Talagrand}, {\em On {R}usso's approximate zero-one law}, Ann. Probab.,
  22 (1994), pp.~1576--1587.

\bibitem{Tas14}
{\sc V.~Tassion}, {\em Planarit\'e et localit\'e en percolation}, PhD thesis,
  ENS Lyon, (2014).

\bibitem{Tas14b}
{\sc V.~Tassion}, {\em Crossing
  probabilities for Voronoi percolation}, Annals of Probability, 44 (2016),
  pp.~3385--3398.
\newblock arXiv:1410.6773.

\bibitem{Why99}
{\sc K.~Whyte}, {\em Amenability, bi-Lipschitz equivalence, and the von Neumann
  conjecture}, Duke Math. J., 99 (1999), pp.~93--112.

\bibitem{yao1977probabilistic}
{\sc A.~C. Yao}, {\em Probabilistic computations: Toward a unified measure of
  complexity}, in Foundations of Computer Science, 1977., 18th Annual Symposium
  on, IEEE, (1977), pp.~222--227.

\end{thebibliography}

\end{document}